\documentclass[11pt,letterpaper]{amsart}
\usepackage{stmaryrd}

\usepackage{tikz}
\usepackage{tikz-cd}

\usepackage{float}
\usepackage{etex}
\usepackage{aliascnt}
\usepackage{bbm} 
\usepackage{pbox}
\usepackage{booktabs}
\usepackage{arydshln}
\usepackage{pinlabel, caption, subcaption}
\usepackage{esvect}
\usepackage{cals}
\usepackage{verbatim}
\usepackage{multicol}
\usepackage[all]{xy}
\usepackage[normalem]{ulem}

\newcommand{\bk}{\mathbbm{k}}

\DeclareMathOperator{\Id}{Id}
\DeclareMathOperator{\Mod}{Mod}

\DeclareMathOperator{\Span}{span}
\DeclareMathOperator{\VB}{VB}

\DeclareMathOperator{\Sp}{Sp}

\newcommand{\VV}{\mathbf{V}}

\newcommand{\HH}{\mathbf{B}}
\newcommand{\HA}{\mathbf{A}}

\newcommand{\cJ}{\mathcal{J}}

\newcommand{\cS}{\mathcal{S}}\newcommand{\cT}{\mathcal{T}}

\providecommand{\genus}{\ensuremath\mr{genus}}

\newcommand{\hooklongrightarrow}{\lhook\joinrel\longrightarrow}

\newcommand{\fL}{\mathfrak{L}}

\newcommand{\gU}{\mathbf{U}}
\newcommand{\gV}{\mathbf{V}}
\newcommand{\gW}{\mathbf{W}}
\newcommand{\gX}{\mathbf{X}}
\newcommand{\gY}{\mathbf{Y}}
\newcommand{\gZ}{\mathbf{Z}}

\newcommand{\F}{\mathbb{F}}

\newcommand{\bQ}{\mathbb{Q}}

\newcommand{\bZ}{\mathbb{Z}}
\newcommand{\Z}{\mathbb{Z}}
\newcommand{\N}{\mathbb{N}_0}
\newcommand{\Q}{\mathbb{Q}}
\newcommand{\C}{\mathbb{C}}

\newcommand{\PP}{\mathfrak{P}}
\newcommand{\LL}{\mathfrak{L}}

\usepackage{lmodern}

\usepackage[dvipsnames,svgnames,x11names]{xcolor}
\usepackage{mathtools,url,graphicx,amssymb,enumerate,nicefrac}
\usepackage[pagebackref,colorlinks,citecolor=blue,linkcolor=blue,urlcolor=blue,filecolor=blue]{hyperref}
\usepackage{microtype}
\usepackage[margin=1.25in]{geometry}
\usepackage{pdflscape}

\numberwithin{theoremcounter}{section}
\newaliascnt{theoremauto}{theoremcounter}

\newaliascnt{Defauto}{theoremcounter}

\newaliascnt{exampleauto}{theoremcounter}

\newaliascnt{lemmaauto}{theoremcounter}

\newaliascnt{propositionauto}{theoremcounter}

\newaliascnt{corollaryauto}{theoremcounter}

\newaliascnt{remarkauto}{theoremcounter}

\newaliascnt{notationauto}{theoremcounter}

\newaliascnt{claimauto}{theoremcounter}

\newaliascnt{warningauto}{theoremcounter}

\newaliascnt{questionauto}{theoremcounter}

\newaliascnt{discussionauto}{theoremcounter}

\newaliascnt{computationauto}{theoremcounter}

\newaliascnt{conjectureauto}{theoremcounter}

\newaliascnt{convauto}{theoremcounter}

\newtheorem{theorem}[theoremauto]{Theorem}
\newtheorem{lemma}[lemmaauto]{Lemma}
\newtheorem{proposition}[propositionauto]{Proposition}
\newtheorem{corollary}[corollaryauto]{Corollary}
\newtheorem*{corollary*}{Corollary}

\theoremstyle{definition}
\newtheorem{definition}[Defauto]{Definition}

\newtheorem{convention}[convauto]{Convention}
\theoremstyle{remark}
\newtheorem{example}[exampleauto]{Example}
\newtheorem{remark}[remarkauto]{Remark}

\newcommand{\mr}[1]{{\rm #1}}

\newcommand{\Stab}{\mr{Stab}}

\DeclareMathOperator*{\sgn}{sgn}

\DeclareMathOperator{\rk}{rk}
\DeclareMathOperator{\GL}{GL}
\DeclareMathOperator{\SL}{SL}
\DeclareMathOperator{\St}{St}

\DeclareMathOperator{\id}{id}

\DeclareMathOperator{\im}{im}

\usepackage{todonotes}

\title[Comodule algebras and the cohomology of symplectic groups]{Comodule algebras and  the  cohomology\\ of symplectic groups}

\author{Jeremy Miller}\thanks{Jeremy Miller was supported by NSF Grant DMS2504473 and a Simons Foundation Travel Support for Mathematicians Grant}
 \email{jeremykmiller@purdue.edu}  
\address{Purdue University \\
	 West Lafayette IN 47907 \\USA}  

\author{Urshita Pal}
\address{2074 East Hall\\
530 Church Street\\
Ann Arbor, MI 48109\\ USA}
\email{urshita@umich.edu}

\author{Peter Patzt}
\thanks{Peter Patzt was supported by NSF Grant DMS2202943 and a Simons Foundation Travel Support for Mathematicians Grant}
\address{University of Oklahoma \\
601 Elm Avenue \\
Norman, OK 73019\\ USA}
\email{ppatzt@ou.edu}

\AtBeginDocument{%
	\def\MR#1{}}

\usepackage{alphabeta}

\begin{document}
	
\begin{abstract} 
We prove that the homology of symplectic groups with Steinberg module coefficients forms a comodule algebra over the homology of general linear groups with Steinberg module coefficients. We prove a freeness-type result for comodule algebras. We deduce new lower bounds for the high-dimensional cohomology of symplectic groups of number rings and their congruence subgroups. 
\end{abstract}

\maketitle

\vspace{-.3in}

\tableofcontents

\section{Introduction}

The goal of this paper is to establish a new structural result for the cohomology of symplectic groups of a Dedekind domain $R$ with coefficients in the symplectic Steinberg module. We apply our results in the case $R$ is a number ring to establish new lower bounds for the cohomology of $\Sp_{2n}(R)$ and its congruence subgroups. Before we state our result, we recall basic definitions and Borel--Serre duality. 

\subsection{Borel--Serre duality}
Let $F$ be the field of fractions of $R$. Let $V$ be an $n$-dimensional $F$-vector space and let $\cT(V)$ denote the geometric realization of the poset of proper nonzero subspaces of $V$ ordered by inclusion. Let $W$ be a $2n$-dimensional symplectic vector space and let $\cT^\omega(W)$ denote the geometric realization of the poset of nonzero isotropic subspace of $W$ ordered by inclusion. The Solomon-Tits theorem implies that both of these spaces are spherical. The Steinberg modules and symplectic Steinberg modules can be defined as the top homology of these simplicial complexes:  $$\St(V) : = \widetilde H_{n-2}(\cT(V)) \text{ and } \St^\omega(W) : = \widetilde H_{n-1}(\cT^\omega(W)).$$ One motivation for considering these representations is their appearance in Borel--Serre duality. Assume $R$ is a number ring with $r_1$ real and $r_2$ many complex conjugate pairs of complex conjugate embeddings. Borel--Serre proved that $\SL_n(R)$ is a virtual duality group of dimension $$\nu^{\SL}_n:= r_1\left(\frac{(n+1)n-2}{2}\right) + r_2(n^2-1) - n + 1 $$ with dualizing module $\St(F^n)$ and  
$\Sp_{2n}(R)$ is a virtual duality group of dimension $$\nu^\omega_{2n}:= r_1 (n(n+1))+r_2(n(2n+1))-n$$ with dualizing module $\St^\omega_{2n}(F^{2n})$, with $F^{2n}$ equipped with the standard symplectic form. 
In particular, $$H^{\nu_{2n}^\omega-i}(\Sp_{2n}(R) ;\Q) \cong H_i(\Sp_{2n}(R) ;\St^\omega(F^{2n}) \otimes \Q).  $$
We call the above groups the \emph{codimension}-$i$ rational cohomology groups since the rational cohomology vanishes above the virtual cohomological dimension. 

Since the homology with Steinberg coefficients determines the untwisted rational cohomology when $R$ is a number ring, there is no loss of information in considering homology with Steinberg coefficients. Moreover, the homology with Steinberg coefficients has additional algebraic structure.

\subsection{Comodule algebras and symplectic groups}
Fix a coefficient ring $\bk$ which we will assume is a field in this section.
Let $$\HH :=\bigoplus_{[P],i} H_i(\GL(P);\St(P \otimes_R F) \otimes \bk)$$ with the sum taken over all isomorphism classes of finitely generated projective $R$-modules and let $$\HA :=\bigoplus_{n,i } H_i(\Sp_{2n}(R);\St^\omega(F^{2n})\otimes \bk).$$ We prove the following.
\begin{theorem}\label{main} Let $R$ be a Dedekind domain. Then $\HA$ is a graded-commutative comodule-algebra over the graded-commutative Hopf algebra $\HH$.
\end{theorem}

That $\HH$ is a graded-commutative Hopf algebra was previously known by work of \cite{AMP} and Brown--Chan--Galatius--Payne \cite{BCGP}. It follows from a theorem of Leray that graded-commutative Hopf algebras over fields of characteristic zero are free as graded-commutative algebras. This allows one to easily produce classes in the rational cohomology of $\SL_n(R)$ for $R$ a number ring. We prove a similar freeness result for comodule algebras \autoref{ComoduleAlgFree}. We use our structural results to deduce the following numerical consequences.

\begin{corollary} \label{corClassNumber} Let $R$ be a number ring with class number $c$. Then $$\dim_{\Q} H^{\nu_{2n}^\omega -\nu^\omega_{2m}}(\Sp_{2n}(R);\Q) \geq (c-1)^{n-m}$$ for all $m \leq n$.


     
     
\end{corollary}

See \autoref{corClassNumberBetter} and \autoref{FrankeRem} for a stronger result along these lines and \autoref{codim96} for an application.
\autoref{corClassNumber} generalizes Br\"uck--Himes \cite[Theorem 1.1]{BH} who proved this result in the case $m=0$. 

\begin{example}
    If $R$ is the ring of integers in $\Q[\sqrt{-23}]$, then $r_1=0$, $r_2=1$, and $c=3$ and hence 
$$\dim_\Q H^{2n^2-2}(\Sp_{2n}(R);\Q) \geq 2^{n-1}$$ for $n \geq 1$. 
\end{example}

\begin{remark}
    Pal--Payne \cite{PalPayne} construct a comodule algebra structure on the Steinberg homology of $\Sp_{2n}(\mathbb{Z})$ using a symplectic analogue of the Waldhausen S-construction, building on arguments from Brown--Chan--Galatius--Payne \cite{BCGP}. They work linear-dually with the corresponding module coalgebra, and identify the action on the coaugmentation class with the weight-zero inclusion in Deligne’s mixed Hodge structure and show that the coproduct is a morphism of mixed Hodge structures. Their results give freeness over the weight zero Hopf algebra and strong lower bounds for Betti numbers and dimensions of weight-graded pieces.
\end{remark}

\subsection{Duoidal comodule algebras and congruence subgroups}
\autoref{main} is the avatar of algebraic structure on the Steinberg modules themselves. The collection of representations $$\{\St(P \otimes_R F)\} $$ has the structure of a graded-commutative duoidal bialgebra and we prove that the collection of representations $$\{\St^\omega(W)\} $$ has the structure of a graded-commutative duoidal comodule algebra over this duoidal bialgebra in \autoref{commutingSp}.

This duoidal perspective is essential for applications to congruence subgroups. For $I$ an ideal in $R$, let $\Sp_{2n}(I)$ denote the kernel of mod-$I$ reduction map $$\Sp_{2n}(R) \to \Sp_{2n}(R/I).$$ When $R=\bZ$, the classifying space $B\Sp_{2n}(j\Z)$ has the homotopy type of the moduli stack of dimension-$n$ principally polarized abelian variety with level-$j$ structure. 
\autoref{commutingSp} has the following numerical consequence.



\begin{corollary} \label{CorCong} Let $p$ be an odd prime and $\bk$ be a field.
     The codimension $i$-cohomology satisfies $$\dim_\bk H^{(n+1)^2 -i}(\Sp_{2n+2}(p\Z);\bk) \geq p^{2n+1}\left(\frac{p-1}{2}\right) \dim_\bk H^{n^2 -i}(\Sp_{2n}(p\Z);\bk).$$ 
     
\end{corollary}

\begin{example}
    Since the $0$th cohomology group of any group is $\bk$, an elementary proof by induction shows that
$$\dim_\bk H^{n^2-m^2}(\Sp_{2n}(p \Z );\bk) \geq p^{n^2-m^2}\left(\frac{p-1}{2}\right)^{n-m}$$ for $n \geq m$.
\end{example}

Capovilla-Searle \cite[Corollary C]{Fabio} obtained the same lower bound for $m=0$.

\subsection{Overview} In \autoref{SecReview}, we recall basic properties of Steinberg modules. In \autoref{DuoidalSection}, we describe the necessary categorical background such as duoidal categories. Then, in \autoref{StDuoidalSec} we prove the symplectic Steinberg modules are duoidal comodule algebras and prove that $\HA$ is a comodule algebra. In \autoref{SecFree}, we prove our freeness results. We apply these results to the cohomology of symplectic groups in \autoref{ApSm} and to their congruence subgroups in \autoref{ApCong}.

\subsection{Acknowledgments} We thank Erik Panzer, Dan Petersen, Andrew Snowden, Robin Sroka, and Jennifer Wilson for helpful conversations. ChatGPT was used for proofreading and reference requests.

\section{Review of Steinberg modules}

In this section, we recall basic facts and constructions related to Steinberg modules.

\label{SecReview}

\subsection{Steinberg modules associated to general linear groups}

We begin by recalling the Steinberg modules associated to the general linear groups and their presentation.

\begin{convention}
  Let $R$ be a Dedekind domain with field of fractions $F$. We will also fix a commutative coefficient ring $\bk$. Also homology groups will be taken with $\bk$-coefficients, all tensor products without subscript will be over $\bk$, etc. 
\end{convention}

\begin{definition} Let $M$ be a projective $R$-module. Let $\cT(M)$ denote the geometric realization of the poset of proper nonzero summands ordered by inclusion.
\end{definition}

\begin{lemma} Let $M$ be a projective module over a Dedekind domain $R$ with field of fractions $F$.
Sending a summand $N$ to $N \otimes_R F$ induces an isomorphism $\cT(M) \to \cT(M \otimes_R F)$.
\end{lemma}

\begin{proof}
    The inverse is given by sending a subspace $V \subset M \otimes_R F$ to $V \cap M$ (compare with Church--Farb--Putman \cite[Lemma 2.2(c)]{CFP}).
\end{proof}

The Solomon--Tits theorem implies $\cT(M)$ is highly connected. Define the rank of $M$ to be $\dim_{F} (M \otimes F)$.

\begin{theorem}[Solomon--Tits] Let $M$ be a rank-$n$ projective module over a Dedekind domain $R$. Then $\cT(M) \simeq \bigvee S^{n-2}$.
\end{theorem}

\begin{definition} Let $M$ be a rank-$n$ projective module over a Dedekind domain $R$. If $n>0$, let $\St(M)=\widetilde H_{n-2}(\cT(M))$. If $n=0$, let $\St(M)=\bk$.
\end{definition}

By convention, $\widetilde H_{-1}(\varnothing) = \bk$. 
The Solomon--Tits theorem also provides generators of the Steinberg modules. Let $\beta=(\vec v_1,\ldots,\vec v_n)$ be an ordered basis of $V=M \otimes_R F$. Let $A(\beta)$ denote the full subcomplex of $\cT(V)$ with vertices subspaces given by the span of proper nonempty subsets of $\{\vec v_1, \dots, \vec v_n\}$. $A(\beta)$ is homeomorphic to $S^{n-2}$ and this sphere comes with an orientation induced by the order on the basis. The apartment class $[\beta]\in \St(V)$ is the image of the fundamental class of $A(\beta)$. The Solomon--Tits theorem states that these are generators of $\St(V) \cong \St(M)$. Lee--Szczarba \cite{LS} determined the relations between these generators.

\begin{theorem}[Lee--Szczarba] \label{GLpres}
Let $V$ be an $n$-dimensional vector space over a field $F$. Then $\St(V)$ has the following presentation. It is  generated by apartment classes $[\vec v_1,\ldots,\vec v_n]$ associated to ordered bases and the relations are as follows:
\begin{enumerate}
   \item Let $u_1,\ldots,u_n \in F^\times$. Then $$[\vec v_1,\ldots,\vec v_n] =[u_1 \vec v_1,\ldots, u_n \vec v_n].$$
     \item Let $\sigma \in S_n$ and let $\sgn\colon S_n \to \{\pm 1\}$ denote the sign function. Then 
     $$[\vec v_1,\ldots,\vec v_n]=\sgn(\sigma) [\vec v_{\sigma(1)},\ldots,\vec v_{\sigma(n)}].$$ 
     \item Let $\vec v_0=\sum_{i=1}^j u_i \vec v_i$ with $u_i \in F^\times$. Then 
     $$\sum_{i=0}^j (-1)^i [\vec v_0,\ldots,\vec v_{i-1},\vec v_{i+1},\ldots,\vec v_j,\ldots,\vec v_n]=0.$$
\end{enumerate}
    
\end{theorem}

\subsection{Steinberg modules associated to symplectic groups}

We now describe the analogous story for symplectic groups. 

\begin{definition}
    Let $W$ be a finite rank free $R$-module. We say that a bi-linear form $\omega$ on $W$ is \emph{symplectic} if there is an ordered basis $(\vec a_1,\vec b_1,\ldots,\vec a_n,\vec b_n)$ with $$\omega(\vec a_i,\vec a_j)=\omega(\vec b_i,\vec b_j)=0 \text{ for all } i \text{ and }j$$ and 
    $$\omega(\vec a_i,\vec b_j)=-\omega(\vec b_i,\vec a_j)=\delta_{i,j} \text{ for all } i \text{ and }j.$$ We say $(\vec a_1,\vec b_1,\ldots,\vec a_n, \vec b_n)$ is a \emph{symplectic basis} and $(W,\omega)$ is a symplectic $R$-module. If $(W,\omega)$ is a symplectic $R$-module, then we say $(\vec x_1,\vec y_1,\ldots,\vec x_n, \vec y_n)$ is a \emph{generalized symplectic basis}  if
    $$\omega(\vec x_i,\vec x_j)=\omega(\vec y_i,\vec y_j)=0 \text{ for all } i \text{ and }j$$ and 
    $$\omega(\vec x_i,\vec y_j)=0 \text{ for all } i \neq j \text{ and } \omega(\vec x_i,\vec y_i)\neq 0.$$ 
\end{definition}

A summand $U \subseteq W$ is called \emph{isotropic} if the form vanishes on $U$ and is called \emph{symplectic} if the form restricts to a symplectic form. Maximal rank isotropic summands are called \emph{Lagrangian}. We call $n$ the \emph{genus} of $W$. If $W$ is a genus $n$ symplectic module over a Dedekind domain $R$ and $J \subset W$ is a rank $m$ isotropic summand, then the symplectic form on $W$ induces a symplectic form on $J^\perp/J$ and $J^\perp/J$ has genus $n-m$. In particular, $J^\perp/J$ is a free $R$-module of rank $2(n-m)$.

\begin{definition}
    Let $(W,\omega)$ be a symplectic module over a Dedekind domain $R$. Let $\cT^\omega(W)$ denote the geometric realization of the poset of nonzero isotropic summands. 
\end{definition}

As before, $\cT^\omega(W) \cong \cT^\omega(W \otimes_R F)$ if $R$ is a Dedekind domain. The Solomon--Tits theorem implies that $\cT^\omega(W)$ is highly connected and provides generators for the top homology.

\begin{theorem}[Solomon--Tits]
Let $(W,\omega)$ be a genus-$n$ symplectic module over a Dedekind domain $R$. Then    $\cT^\omega(W) \simeq \bigvee S^{n-1}$.
\end{theorem}

Let $(W,\omega)$ be a genus-$n$ sympletic module over a field $F$. Given an ordered generalized symplectic basis $\beta=(\vec a_1,\vec b_1,\ldots,\vec a_n,\vec b_n)$ of $(W,\omega)$, let $B(\beta)$ denote the full subcomplex of $\cT^\omega(W)$ on vertices corresponding to nonzero isotropic summands that are equal to the span of a subset of $\beta$. Note, unlike in the case of general linear groups, we do not take all proper nonempty subsets, as not all such subsets span isotropic subspaces. We have that $B(\beta) \cong S^{n-1}$ with a natural orientation coming from the order. Let $[\beta]$ denote the image of fundamental class of $B(\beta)$ in $\St^\omega(W)$.

Let $S_{2n}^\omega \leq S_{2n}$ be the subgroup of permutations such that induce permutations on the following set of subsets $$\{\{1,2 \},\{3,4 \},\dots,\{2n-1,2n \} \}.$$ Abstractly,  $S_{2n}^\omega \cong S_2 \wr S_n$. Let $$\chi\colon S_{2n}^\omega \to \{\pm 1\}  $$ be the product of the sign representations of each $S_2$ and the sign representation of $S_n$.

The following presentation can be extracted from \cite[Corollary 1.4, Lemma 4.3, Proof of Theorem 1.5]{Pal}; see also T\'oth \cite[Theorem 2]{Toth}.

\begin{theorem}[Pal, T\'oth] \label{SympPres}
    Let $(W,\omega)$ be a genus $n$ symplectic vector space over a field $F$. Then $\St^\omega(W)$ has the following presentation. The generators are given by
    symplectic apartment classes $[\vec a_1,\vec b_1,\ldots,\vec a_n,\vec b_n]$ associated to ordered generalized symplectic  bases and the relations are as follows: 
\begin{enumerate}
   \item Let $u_1,\ldots,u_{2n} \in F^\times$. Then $$[\vec a_1,\vec b_1, \ldots,\vec a_n,\vec b_n] =[u_1 \vec a_1,u_2 \vec b_1, \ldots,u_{2n-1} \vec a_n, u_{2n} \vec b_n].$$
     \item Let $\sigma \in S^{\omega}_{2n}$.  Then
     $$[\vec v_1, \ldots,\vec v_{2n}]=\chi(\sigma)[\vec v_{\sigma(1)}, \ldots,\vec v_{\sigma(2n)}].$$ 
     \item    
     $ [\vec a_1,\vec b_1, \vec a_2, \vec b_2,\ldots,\vec a_n,\vec b_n  ]=$ $$[\vec a_1+\vec b_1,\vec b_1, \vec a_2,\vec b_2,\ldots,\vec a_n,\vec b_n  ]+[\vec a_1,\vec a_1+\vec b_1, \vec a_2, \vec b_2,\ldots,\vec a_n,\vec b_n  ].$$
     \item $ [\vec a_1,\vec b_1, \vec a_2, \vec b_2,\ldots,\vec a_n,\vec b_n  ]=$  $$ [\vec a_1,\vec b_1-\vec b_2, \vec a_1+\vec a_2, \vec b_2,\ldots,\vec a_n,\vec b_n  ]+[\vec b_1-\vec b_2,\vec a_2, \vec a_1+\vec a_2, \vec b_1,\ldots,\vec a_n,\vec b_n  ].$$
\end{enumerate}
\end{theorem}


\section{Algebraic structure}
In this section, we describe monoidal structures on categories of representations of general linear groups and symplectic groups. 

\label{DuoidalSection}
\subsection{Monoidal structures}

\begin{definition}
    For $R$ a ring, let $\VB(R)$ denote the groupoid of finitely generated projective $R$-modules and isomorphisms. Let $\VB^\omega(R)$ denote the groupoid of finite genus symplectic $R$-modules and form-preserving isomorphisms.
\end{definition}

\begin{definition}
    Let $R$ and $\bk$ be rings. A $\VB(R)$-module over $\bk$ is a functor $\gX\colon \VB(R) \to \Mod_\bk$. Let $\Mod_{\VB(R)}(\bk)$ denote the category of $\VB(R)$-modules over $\bk$.  
    Similarly define  $\VB(R)$-chain complexes and the corresponding symplectic versions. 
\end{definition}

\begin{definition} Given two $\VB(R)$-modules $\gX$ and $\gY$, let $\gX \otimes_\LL \gY$ be the $\VB(R)$-module with value on a projective $R$-module $M$ is given by $$(\gX \otimes_\LL \gY)(M) =\bigoplus_{A \oplus B=M} \gX(A) \otimes \gY(B).$$ 
Let  $\gX \otimes_\PP \gY$ be the $\VB(R)$-module with value on a projective $R$-module $M$ is given by $$(\gX \otimes_\PP \gY)(M) =\bigoplus_{A \subseteq M} \gX(A) \otimes \gY(M/A).$$ Here $A$ runs over all summands of $M$. 
\end{definition}

The products $\otimes_\LL$ and $\otimes_\PP$ endow $\Mod_{\VB(R)}(\bk)$ with monoidal structures. The unit of $\otimes_\LL$ and $\otimes_\PP$ are given by the $\VB(R)$-module $$M \mapsto \begin{cases}
\bk, & \text{if } M\cong 0,\\
0, & \text{otherwise.}
\end{cases} $$ While $\otimes_\LL$ is symmetric monoidal, $\otimes_\PP$ is not. These two structures are compatible in the sense that they form a duoidal category. This will be discussed in \autoref{duoidal cat}. Part of this data is an interchange natural transformation.

\begin{definition}
    Let $\gX,\gY,\gZ,\gW$ be $\VB(R)$-modules. Let 
    $$\zeta_{\gX,\gY,\gZ,\gW} \colon (\gX \otimes_\PP \gY) \otimes_\LL (\gZ \otimes_\PP \gW) \to (\gX \otimes_\LL \gZ) \otimes_\PP (\gY \otimes_\LL \gW)$$ 
    be the following map. For $M$ a finitely generated projective $R$-module, 
    \begin{align*}
        &\Big((\gX \otimes_\PP \gY) \otimes_\LL (\gZ \otimes_\PP \gW)\Big)(M)\\
        =&\bigoplus_{A \oplus B =M} (\gX \otimes_\PP \gY)(A) \otimes (\gZ \otimes_\PP \gW)(B)\\
        =&\bigoplus_{\substack{A \oplus B =M\\  A' \subseteq A,\, B' \subseteq B}} \gX(A') \otimes \gY(A/A') \otimes \gZ(B') \otimes \gW(B/B'). 
    \end{align*}
    Similarly, 
    \begin{align*}
    &\big((\gX \otimes_\LL \gZ) \otimes_\PP (\gY \otimes_\LL \gW) \Big)(M)\\
    =&\bigoplus_{M' \subseteq M} (\gX \otimes_\LL \gZ)(M') \otimes (\gY \otimes_\LL \gW)(M/M')\\
    =&\bigoplus_{\substack{ M' \subseteq M\\C \oplus D=M', \, C' \oplus D'=M/M'} } \gX(C) \otimes \gZ(D) \otimes \gY(C') \otimes \gW(D').
    \end{align*}
    A collection $A,B,A',B'$  in the first indexing set determines a collection $C,D,C',D',M'$ in the second indexing set as follows.
    Let  $M'=A' \oplus B', C=A'$,  
    $D=B'$, $ 
    C'=A/A' $, $D'=B/B'$ and identify $M/(A' \oplus B)'
    $ with $(A/A') \oplus (B/B')$. The identity maps $A' \to C$ and $B' \to D$, the inverses of the natural isomorphism $ 
    C' \cong A/A'$ and $D' \cong B/B'$, and the symmetry isomorphism  $\gY(C') \otimes \gZ(D) \cong \gZ(D) \otimes \gY(C')$    
    induce $\zeta_{\gX,\gY,\gZ,\gW}$. 
\end{definition}

Note that the symmetry isomorphism $\gY(C') \otimes \gZ(D) \cong \gZ(D) \otimes \gY(C')$ includes a sign of  $(-1)^{\rk(C')\rk(D)}$.

There are symplectic versions of $\otimes_\LL$ and $\otimes_\PP$. Recall that for $(W,\omega)$ a symplectic module and $M \subset W$, the \emph{symplectic complement} 
$$M^\perp=\{w \in W \,|\, \omega(w,m) =0\text{ for all } m \in M \}.$$ 
If $M$ is symplectic, then $M^\perp$ is symplectic and $W =M \oplus M^\perp$. If $R$ is a Dedekind domain, then the genus of $W$ is the sum of the genus of $M$ and the genus of  $M^\perp$. If $M$ is isotropic $M \subset M^\perp$ and $M^\perp/M$ has a natural symplectic structure. If $R$ is a Dedekind domain, then the genus of $M^\perp /M$ is the genus of $W$ minus rank of $M$.

\begin{definition}
    Let $\gU$ and $\gV$ be $\VB^\omega(R)$-modules and let $\gX$ be a $\VB(R)$-module. Let $\gU \otimes_\LL^\omega \gV$ be the $\VB^\omega(R)$-module with value on a symplectic module $W$ given by $$(\gU \otimes_\LL^\omega \gV)(W)=\bigoplus_{A \oplus B=W}\gU(A) \otimes \gV(B).$$ Here $A$ and $B$ are symplectic submodules of $W$ and $B=A^\perp$. $\Mod_{\VB^\omega(R)}(\bk)$ is monoidal category with $\otimes_\LL^\omega$.

    Let $\gX \otimes_\PP^\omega \gV$ be the $\VB(R)^\omega$-module with value on a symplectic module $W$ given by $$(\gX \otimes_\PP^\omega \gV)(W)=\bigoplus_{J \subset W}\gX(J) \otimes \gV(J^\perp/J).$$ Here $J$ is an isotropic summand. The category and operation $(\Mod_{\VB^\omega}(R)(\bk),\otimes_\LL^\omega)$ has the structure of a monoidal  module category over $\Mod_{\VB(R)}(\bk)$ (sometimes called a $\Mod_{\VB(R)}(\bk)$-actegory).
\end{definition}

We now describe a symplectic analogue of the interchange.

\begin{definition}
    Let $\gX$ and $\gZ$ be $\VB(R)$-modules and let $\gY$ and $\gW$ be $\VB^\omega(R)$-modules. Let 
    $$\xi_{\gX,\gY,\gZ,\gW} \colon (\gX \otimes_\PP^\omega \gY) \otimes_\LL^\omega (\gZ \otimes_\PP^\omega \gW) \to (\gX \otimes_\LL \gZ) \otimes_\PP^\omega (\gY \otimes_\LL^\omega \gW)$$ 
    be the following map. For $W$ a finitely generated symplectic $R$-module, 
    $$\Big((\gX \otimes^\omega_\PP \gY) \otimes_\LL^\omega (\gZ \otimes_\PP^\omega \gW)\Big)(W)=$$
    $$\bigoplus_{U \oplus V =W} (\gX \otimes^\omega_\PP \gY)(U) \otimes (\gZ \otimes_\PP^\omega \gW)(V) = $$
    $$\bigoplus_{\substack{U \oplus V =W\\  J \subseteq U,\, K \subseteq V}} \gX(J) \otimes \gY(J^{\perp U}/J) \otimes \gZ(K) \otimes \gW(K^{\perp V}/K).  $$
    Here $U,V$ are symplectic, $J,K$ are isotropic, and $J^{\perp U}$ and $K^{\perp V}$ are the complements of $J$ in $U$ and $K$ in $V$, respectively.  
    Similarly, 
    $$\big((\gX \otimes_\LL \gZ) \otimes_\PP^\omega (\gY \otimes_\LL^\omega \gW) \Big)(W)=$$ 
    $$\bigoplus_{Q \subset W} (\gX \otimes_\LL \gZ)(Q) \otimes (\gY \otimes_\LL^\omega \gW)(Q^\perp/Q)=$$ 
    $$\bigoplus_{\substack{ Q \subset  W\\C \oplus D=Q, \, C' \oplus D'=Q^\perp/Q} } \gX(C) \otimes \gZ(D) \otimes \gY(C') \otimes \gW(D').$$ 
    Here $Q$ is isotropic and $C',D'$ are symplectic.
A collection $U,V,J,K$ in the first indexing set determines a collection $Q,C,D,C',D'$ in the second indexing set as follows. Let $Q=J \oplus K$, $C=J$, $D=K$, $C'=J^{\perp U}/J$ and $D'=K^{\perp V}/K$. Note that $J^{\perp U}/J \oplus K^{\perp V}/K = Q^\perp/Q$, where $Q^\perp$ is the complement of $Q$ in $W$. The identity maps $J \to C$, $K \to D$, the projection isomorphisms $J^{\perp U} \to C'$, $K^{\perp V}/K \to D'$, and the symmetry  $\gY(C') \otimes \gZ(D) \to \gZ(D) \otimes \gY(C')$ induce $\xi_{\gX,\gY,\gZ,\gW}$. 
\end{definition}

Again, the symmetry isomorphism $\gY(C') \otimes \gZ(D) \to \gZ(D) \otimes \gY(C')$ includes a sign of  $(-1)^{\genus(C')\rk(D)}$.

\subsection{Duoidal compatibility} \label{duoidal cat}

Duoidal categories provide a framework for describing the relationship between $\PP$ and $\LL$.

\begin{definition}
A \emph{duoidal category} is a category $\mathcal{C}$ equipped with two monoidal structures $(\mathcal{C}, \diamond, I)$ and $(\mathcal{C}, \star, J)$, such that the functor $\star \colon \mathcal{C} \times \mathcal{C} \to \mathcal{C}$ and the object $J \colon \mathbf{1} \to \mathcal{C}$ are lax monoidal functors with respect to the monoidal structure $(\mathcal{C}, \diamond, I)$.
\end{definition}

The data of a duoidal category is the following:
\begin{enumerate}
    \item An \emph{interchange} natural transformation:
    \[ \zeta_{X,Y,Z,W} \colon (X \star Y) \diamond (Z \star W) \longrightarrow (X \diamond Z) \star (Y \diamond W) \]
    
    \item A morphism for the unit of $\diamond$:
    \[ \Delta_I \colon I \longrightarrow I \star I \]
    
    \item A morphism for the unit of $\star$:
    \[ \mu_J \colon J \diamond J \longrightarrow J \]
    
    \item A morphism relating the two units:
    \[ \delta \colon I \longrightarrow J \]
\end{enumerate}

These morphisms satisfy compatibility axioms. We will often omit units and counits from the notation.

See Kupers--Rudenko--Sierra \cite[Section 2.2]{KRS1} for a proof that 
$(\Mod_{\VB(R)}(\bk),\otimes_\PP,\otimes_\LL)$ is a duoidal category in the case $R$ is a field. A similar argument works for Dedekind domains.
\begin{definition}
    A duoidal bialgebra in a duoidal category $(\mathcal{C}, \star, J, \diamond, I,\zeta)$ is an object $B$, a monoid structure $(B,\mu,\nu)$ on $B$ with respect to $\diamond$, and a comonoid structure $(B,\Delta,\epsilon)$ with respect to $\star$ such that the following diagrams commutes:
\begin{enumerate}
    \item   
        \begin{tikzcd}
        B \diamond B \arrow[d, "\mu"] \arrow[r, "\Delta \diamond \Delta"]&  (B \star B) \diamond (B \star B) \arrow[d, "(\mu \star \mu) \circ \zeta"] \\ 
            B  \arrow[r, "\Delta"] & B \star B,
        \end{tikzcd}
    \item  \begin{tikzcd}
        I \arrow[d, "\nu "]  & \arrow[l, "\mu_I"] I \diamond I  \arrow[d, "\nu \diamond \nu"]    \\ 
            B   & B \diamond B,   \arrow[l, "\mu"]
        \end{tikzcd}
    \item \begin{tikzcd}
        J    \arrow[r, "\Delta_J "] & J \star J    \\ 
            B \arrow[u, "\epsilon "]   \arrow[r, "\Delta"] & \arrow[u, "\epsilon \star \epsilon"] B \star B,
        \end{tikzcd}
    \item \begin{tikzcd}
        I   \arrow[d, "\delta "] \arrow[dr, " \nu"]   \\ 
            J   & B.  \arrow[l, "\epsilon"]
        \end{tikzcd}
\end{enumerate}

\end{definition}

\begin{definition}
A \emph{duoidal module category} over a duoidal category  $(\mathcal{C},  \star^c, J^c , \diamond^c, I^c)$ is a monoidal category $(\mathcal{D},\diamond^d, I^d)$, a left $(\mathcal{C},\star^c,J^c)$-module structure $\star^d \colon \mathcal{C} \times \mathcal{D} \to \mathcal{D}$ 
such that the functor $\star^d \colon \mathcal{C} \times \mathcal{D} \to \mathcal{D}$ 
is lax monoidal with respect to the monoidal structure coming from $(\mathcal{C}, \diamond^c, I^c)$ and $(\mathcal{D}, \diamond^d, I^d)$.

Explicitly, this means there exist the following structural morphisms:
\begin{enumerate}
    \item An \emph{interchange} natural transformation:
    \[ \xi_{X,Y,Z,W} \colon (X \star^d Y) \diamond^d (Z \star^d W) \longrightarrow (X \diamond^c Z) \star^d (Y \diamond^d W) \]
    
    \item A morphism for the unit of $\diamond^d$:
    \[ \Delta_I \colon I^d \longrightarrow I^c \star^d I^d \]
    
    
\end{enumerate}
\noindent These structural morphisms must satisfy a set of coherence axioms ensuring compatibility with the associativity and unitality constraints of the monoidal and module structures.
\end{definition}

A similar argument to that of Kupers--Rudenko--Sierra \cite{KRS1} gives the following.

\begin{proposition} \label{commutingGL}
For $R$ a Dedekind domain, $(\Mod_{\VB^\omega(R)}(\bk),\PP,\LL)$ is a duoidal module category over  $(\Mod_{\VB(R)}(\bk),\PP,\LL)$.
\end{proposition}

\begin{definition}  \label{duoidalcomodulealgebra}
Let  $(\mathcal{D},\diamond^d, I^d,\xi)$ be a  duoidal module category over a duoidal category $(\mathcal{C}, \diamond^c, I^c, \star^c, J^c)$.   Let $(B,\mu_B,\Delta_B,\epsilon_B,\nu_B)$ be a duoidal bialgebra in $(\mathcal{C}, \diamond^c, I^c, \star^c, J^c)$. Let $(A,\mu_A,\nu_A)$ be a monoid object in $(\mathcal{D},\diamond^d, I^d)$ with a comodule structure map $\Delta_A\colon A \to B \star^d A$. Then $A$ is called a duoidal comodule algebra if the following diagrams commute:
\begin{enumerate}
    \item    \label{duoidalcomodulealgebra1}
        \begin{tikzcd}
        A \diamond^d A \arrow[d, "\mu_A"] \arrow[r, "\Delta_A \diamond \Delta_A"]&  (B \star_d A) \diamond^d (B \star_d A) \arrow[d, "(\mu_B \star^d \mu_A) \circ \xi"] \\ 
            A  \arrow[r, "\Delta_A"] & B \star^d A,
        \end{tikzcd}
    \item  \begin{tikzcd} 
        I^d \arrow[d, "\eta_A "]\arrow[r, "\Delta^d_{I} "] & I^d \star^d I^d \arrow[d, "\nu_B \diamond^d \nu_A"]      \\ 
            A  \arrow[r, "\Delta_A"]   & B \star^d A,  
        \end{tikzcd}
\end{enumerate}

\end{definition}

\section{Algebraic Structure of Steinberg modules}
\label{StDuoidalSec}

In this section, we will prove that symplectic Steinberg modules are a duoidal comodule algebra and their homology is a comodule algebra. 

\subsection{Structure maps}

We first recall the product and coproduct structure on general linear group Steinberg modules. Then we extend these constructions to symplectic Steinberg modulues. 
\begin{definition}
    Let $\St$ denote the $\VB(R)$-module with value on $M$ given by $\St(M)$.     Let $\St^\omega$ denote the $\VB^\omega(R)$-module with value on $W$ given by $\St^\omega(W)$.
\end{definition}

For $\VB$-modules $\gX,\gY,\gZ$, note that a morphism $$f\colon \gX \otimes_\LL \gY \to \gZ $$ is determined by the data of $\GL(M) \times \GL(N)$-equivariant linear maps $$f_{M,N}\colon \gX(M) \otimes \gY(N) \to \gZ(M \oplus N).$$

\begin{definition}
    Let $\mu\colon \St \otimes_\LL \St \to \St$ have structure maps 
    \begin{align*}
        \mu_{M,N}\colon \St(M) \otimes \St(N) &\longrightarrow \St(M \oplus N)\\
        [\vec v_1,\ldots, \vec v_p ] \otimes [\vec w_1,\ldots, \vec w_q ] &\longmapsto [\vec v_1,\ldots, \vec v_p, \vec w_1,\ldots, \vec w_q ]. 
    \end{align*}    
    Here $\vec v_1,\ldots,\vec v_p$ is a basis of $M \otimes_R F$ and $\vec w_1,\ldots,\vec w_q$ is a basis of $N \otimes_R F$. 
\end{definition}

See \cite{MNP}, \cite{AMP}, and Kupers--Rudenko--Sierra \cite{KRS1} for a proof that $\mu$ is well defined, associative and graded-commutative. Given a $\bk$-module $T$, let $T[0]$ denote the $\VB(R)$-module with value $T$ on rank-$0$ $R$-modules and value $0$ otherwise. Note that $\bk[0]$ is the unit of $\Mod_{\VB(R)}(\bk)$. The unit of $\St$ is the map $\bk \to \St$ induced by the isomorphism $\bk \to \St(0)$.

For $\VB$-modules $\gX,\gY,\gZ$, note that a morphism $$f\colon  \gX \to \gY \otimes_\PP \gZ  $$ is the same data as linear maps $$f_{M' \subseteq M}\colon  \gX(M) \to \gY(M') \otimes \gZ(M /M')$$ which are equivariant with respect to the parabolic subgroup that is the stabilizer of $M'$ in $\GL(M)$ and which are zero for all but finitely many summands $M'$ of $M$.

\begin{definition}
    Let $\Delta\colon  \St \to \St \otimes_\PP \St $ have structure maps $$\Delta_{M' \subseteq M}\colon  \St(M) \to \St(M') \otimes \St(M/M') $$ given by the following formula. Let $\vec v_1,\ldots,\vec v_n$ be a basis of $M \otimes_R F$. If there is a shuffle $\sigma \in S_n$ and $m \leq n$ with $$M'=\Span(\vec v_{\sigma(1)},\ldots,\vec v_{\sigma(m)}) \cap M,$$ then let  $$\Delta_{M' \subseteq M} \left( [\vec v_1,\ldots,\vec v_n] \right) = \sgn(\sigma)[\vec v_{\sigma(1)},\ldots,\vec v_{\sigma(m)}] \otimes [\vec v_{\sigma(m+1)},\ldots,\vec v_{\sigma(n)}]. $$ Here we view $\vec v_{\sigma(m+1),\ldots,\vec v_{\sigma(n)}}$ as vectors in $(M/M') \otimes_R F$. Otherwise, let $$\Delta_{M' \subseteq M} \left([v_1,\ldots,v_n]\right) = 0. $$
    
\end{definition}

See \cite{AMP} and Kupers--Rudenko--Sierra \cite{KRS1} for a proof that $\Delta$ is well defined and coassociative. We now describe symplectic versions. 

For $\VB^\omega$-modules $\gX,\gY,\gZ$, note that a morphism $$f\colon \gX \otimes_\LL^\omega \gY \to \gZ $$ is the same data as $\Sp(W) \times \Sp(V)$-equivariant linear maps $$f_{W,V}\colon \gX(W) \otimes \gY(V) \to \gZ(W \oplus V).$$

\begin{definition}
    Let $\mu\colon  \St^\omega \otimes^\omega_\LL \St^\omega \to \St^\omega$ have structure maps $$\mu_{W,V}\colon  \St^\omega(W) \otimes \St^\omega(V) \to \St^\omega(W \oplus V) $$  $$[\vec a_1,\vec b_1\ldots, \vec a_p, \vec b_p ] \otimes [\vec c_1,\vec d_1\ldots, \vec c_q, \vec d_q ] \mapsto [\vec a_1,\vec b_1\ldots, \vec a_p, \vec b_p , \vec c_1,\vec d_1\ldots, \vec c_q, \vec d_q ]. $$ Here $\vec a_1,\vec b_1\ldots, \vec a_p, \vec b_p $ is a generalized symplectic basis of $W \otimes_R F$ and 
    $\vec c_1,\vec d_1\ldots, \vec c_q, \vec d_q $ is a generalized symplectic basis of $V \otimes_R F$.
\end{definition}

We now check that this is well defined. 

\begin{proposition}
The map $\mu_{W,V}\colon  \St^\omega(W) \otimes \St^\omega(V) \to \St^\omega(W \oplus V)$ is well defined.
\end{proposition}
\begin{proof}
Let $\vec a_1,\vec b_1\ldots, \vec a_p, \vec b_p $ be a symplectic basis of $W \otimes_R F$ and 
    $\vec c_1,\vec d_1\ldots, \vec c_q, \vec d_q $ be a symplectic basis of $V \otimes_R F$. If you apply one of the relations in \autoref{SympPres} to $[\vec a_1,\vec b_1\ldots, \vec a_p, \vec b_p]$ or $[\vec c_1,\vec d_1\ldots, \vec c_q, \vec d_q ]$, then you can undo that by performing the corresponding relation (and possibly also relation (2)) to $[\vec a_1,\vec b_1\ldots, \vec a_p, \vec b_p,\vec c_1,\vec d_1\ldots, \vec c_q, \vec d_q ]$. 
\end{proof}
\begin{proposition}
$(\St^{\omega},\mu)$ is a unital graded-commutative monoid object in $(\Mod_{\VB^\omega(R)},\otimes_\LL^\omega)$.    
\end{proposition}

\begin{proof} The image of the unit map is the empty apartment class and concatenation with the empty apartment class is the identity operation. Thus $\St^\omega$ is unital.
    Associativity follows from the fact that 
    \begin{align*}
        &\mu\Big(\mu([\vec a_1,\ldots,\vec a_{2n}],[\vec b_1,\ldots,\vec b_{2m}]),[\vec c_1,\ldots,\vec c_{2k}]\Big)\\
        =&[\vec a_1,\ldots,\vec a_{2n},\vec b_1,\ldots,\vec b_{2m},\vec c_1,\ldots,\vec c_{2k}]\\
        =& \mu \Big([\vec a_1,\ldots,\vec a_{2n}],\mu([\vec b_1,\ldots,\vec b_{2m}],[\vec c_1,\ldots,\vec c_{2k}])\Big).
    \end{align*}    
    Graded commutativity follows from the fact that $$\mu([\vec a_1,\ldots,\vec a_{2n}],[\vec b_1,\ldots,\vec b_{2m}])=(-1)^{nm}\mu([\vec b_1,\ldots,\vec b_{2m}],[\vec a_1,\ldots,\vec a_{2n}]).$$
\end{proof}

For $\VB^\omega$-modules $\gX,\gZ$, and a $\VB$-module $\gY$, note that a morphism $$f\colon  \gX \to \gY \otimes_\PP^\omega \gZ  $$ is determined by linear maps $$f_{J \subseteq W}\colon  \gX(W) \to \gY(J) \otimes \gZ(J^\perp /J)$$ which are equivariant with respect to the parabolic subgroup that is the stabilizer of $J$ in $\Sp(W)$ and which are zero for all but finitely many isotropic summands $J$ of $W$.


\begin{definition}
    Let $\Delta\colon  \St^\omega \to \St \otimes_\PP^\omega \St^\omega $ have structure maps associated to an isotropic $J$ $$\Delta_{J \subset W}\colon  \St^\omega(W) \to \St(J) \otimes \St^\omega(J^\perp/J) $$ given by the following formula. Let $(\vec v_1,\dots,\vec v_{2n})$ be an ordered symplectic basis of $W \otimes_R F$. If there is  no  $\sigma \in S_{2n}^\omega$ and $m \leq n$ with $$J=\Span(\vec v_{\sigma(1)},\vec v_{\sigma(3)},\dots,\vec v_{\sigma(2m-1)}) \cap W,$$ then let $$\Delta_{J \subset W} \left([\vec v_1,\dots,\vec v_{2n}]\right) = 0. $$ Otherwise, fix such a $\sigma$.
    For $i > 2m$, view $\vec v_{\sigma(i)}$ as an element of the quotient space $J^\perp/J$. Let 
    $$\Delta_{J \subset W} \left( [\vec v_1,\dots,\vec v_{2n}] \right) = \chi(\sigma)[\vec v_{\sigma(1)},v_{\sigma(3)},\dots, \vec v_{\sigma(2m-1)} ] \otimes [\vec v_{\sigma(2m+1)},\vec v_{\sigma(2m+2)}, \dots, \vec v_{\sigma(2n-1)}, \vec v_{\sigma(2n)} ]. $$    
    Note that this equation is independent of which $\sigma\in S_{2n}^\omega$ satisfying
    $$J=\Span(\vec v_{\sigma(1)},\vec v_{\sigma(3)},\dots,\vec v_{\sigma(2m-1)}) \cap W$$
    we use.
\end{definition}

We now check that $\Delta$ is well defined. 

\begin{proposition}
    The map $\Delta_{J \subset W}\colon  \St^\omega(W) \to \St(J) \otimes \St^\omega(J^\perp/J)$ is well defined.
\end{proposition}

\begin{proof} Let $\beta=(\vec a_1,\vec b_1,\ldots,\vec a_n,\vec b_n)$ be an ordered generalized symplectic basis of $W$. It is clear that if we apply relations $(1)$ and $(2)$ from \autoref{SympPres} to $\beta$, then the output of $\Delta_{J \subset W}$ is changed by apply relations $(1)$ and $(2)$ from \autoref{GLpres} and \autoref{SympPres}. 

Now apply relation $(3)$ from \autoref{SympPres} to $\beta$. The case $J=0$ is easy so assume $J \neq 0$. By applying relation $(2)$ from \autoref{SympPres}, we may assume we are in of the following cases:
\begin{enumerate}
    \item [i.]  $J \otimes_R F=\Span(\vec a_1,\ldots,\vec a_m)$,
    \item [ii.] $J \otimes_R F=\Span(\vec a_1+\vec b_1,\vec a_2\ldots,\vec a_m)$,
 \item [iii.]  $J \otimes_R F=\Span(\vec a_2,\ldots,\vec a_{m+1})$,
    \item [iv.] $J \otimes_R F=\Span(\vec a_2+\vec b_2,\vec a_3\ldots,\vec a_{m+1})$,
    \item [v.] $J \otimes_R F$ is not the span of a subset of $\beta \cup \{\vec a_i + \vec b_i\}$.
\end{enumerate} In case (i), 
\begin{gather*}\Delta_{J \subset W} ([\vec a_1,\vec b_1, \vec a_2, \vec b_2,\ldots,\vec a_n,\vec b_n  ])=[\vec a_1,\ldots, \vec a_m] \otimes [\vec a_{m+1}, \vec b_{m+1},\ldots,\vec a_n,\vec b_n  ],\\
\Delta_{J \subset W} ([\vec a_1+\vec b_1,\vec b_1, \vec a_2,\vec b_2,\ldots,\vec a_n,\vec b_n  ])=0, \text{ and}\\
\Delta_{J \subset W} ([\vec a_1,\vec a_1+\vec b_1, \vec a_2, \vec b_2,\ldots,\vec a_n,\vec b_n  ])=[\vec a_1,\ldots, \vec a_m] \otimes [\vec a_{m+1}, \vec b_{m+1},\ldots,\vec a_n,\vec b_n  ].
\end{gather*}
This means that applying relation (3) does not change the output in case (i). A similar argument works in case (ii). Now consider case (iii). In case (iii), $$\Delta_{J \subset W} ([\vec a_1,\vec b_1, \vec a_2, \vec b_2,\ldots,\vec a_n,\vec b_n  ])=[\vec a_2,\ldots, \vec a_{m+1}] \otimes [\vec a_1, \vec b_1,\vec a_{m+2}, \vec b_{m+2},\ldots,\vec a_n,\vec b_n  ],$$
$$\Delta_{J \subset W} ([\vec a_1+\vec b_1,\vec b_1, \vec a_2,\vec b_2,\ldots,\vec a_n,\vec b_n  ])=[\vec a_2,\ldots, \vec a_{m+1}] \otimes [\vec a_1+\vec b_1, \vec b_1,\vec a_{m+2}, \vec b_{m+2},\ldots,\vec a_n,\vec b_n  ], \text{ and}$$ 
$$\Delta_{J \subset W} ([\vec a_1,\vec a_1+\vec b_1, \vec a_2, \vec b_2,\ldots,\vec a_n,\vec b_n  ])=[\vec a_2,\ldots, \vec a_{m+1}] \otimes [\vec a_1, \vec a_1+\vec b_1,\vec a_{m+2}, \vec b_{m+2},\ldots,\vec a_n,\vec b_n  ].$$ By relation (3) from \autoref{SympPres}, applying relation (3) from \autoref{SympPres} does not change the output in case (iii). A similar argument works in case (iv).
In case (v), all terms map to $0$. 

Now consider relation (4) from \autoref{SympPres}. The only interesting case to consider is when $J \otimes_R F=\Span(a_2,\ldots,a_{m+1})$ as all other cases can be dealt with by arguments similar to uses previously used in this proof.  We have 
\begin{align*}
    &\Delta_{J \subset W} [\vec a_1,\vec b_1, \vec a_2, \vec b_2,\ldots,\vec a_n,\vec b_n  ]\\
    =&[\vec a_2,\ldots,\vec a_{m+1}] \otimes [\vec a_1,\vec b_1,\vec a_{m+2}, \vec b_{m+2},\ldots,\vec a_n,\vec b_n],\\ 
    &\Delta_{J \subset W} [\vec a_1,\vec b_1-\vec b_2, \vec a_1 + \vec a_2, \vec b_2,\vec a_3,\vec b_3,\ldots,\vec a_n,\vec b_n  ]\\
    =&0, \text{ and}\\
    &\Delta_{J \subset W}([\vec b_1-\vec b_2,\vec a_2, \vec a_1+\vec a_2, \vec b_1,\ldots,\vec a_n,\vec b_n  ])\\
    = &[\vec a_2,\ldots,\vec a_{m+1}] \otimes [\vec a_1+\vec a_2,\vec b_1,\vec a_{m+2}, \vec b_{m+2},\ldots,\vec a_n,\vec b_n]\\
    =&[\vec a_2,\ldots,\vec a_{m+1}] \otimes [\vec a_1,\vec b_1,\vec a_{m+2}, \vec b_{m+2},\ldots,\vec a_n,\vec b_n].
\end{align*}
Here we use that the symplectic apartment in the output of $\Delta_{J \subset W}$ lives in the symplectic Steinberg module of the quotient space $J^\perp/J$. This checks that applying relation (4) from \autoref{SympPres} does not change the output.
\end{proof}

\begin{proposition}
    $(\St^\omega,\Delta)$ is counital  coassociative comodule over $(\St,\Delta)$. 
\end{proposition}
\begin{proof}
Using the fact that the unit maps to the empty apartment class, it is not hard to see that $(\St^\omega,\Delta)$ is counital over $(\St,\Delta)$.

Let $W$ be a symplectic $R$-module and let  $v_1,v_2, \ldots, v_{2n-1},v_{2n}$ be a symplectic basis of $W \otimes_R F$. To check that $\Delta$ is coassociative, We must check that \begin{multline*}
(\Id \otimes \Delta)\Big(\Delta([\vec v_1,\vec v_2, \ldots, \vec v_{2n-1},\vec v_{2n}])\Big)=(\Delta \otimes \Id)\Big(\Delta([\vec v_1,\vec v_2, \ldots, \vec v_{2n-1},\vec v_{2n}])\Big) \\\in \bigoplus_{M \subseteq J \subset W} \St(M) \otimes \St(J/M) \otimes \St^\omega(J^\perp/J).
\end{multline*}
Here $J$ is an isotropic summand of $W$ and $M$ is a summand of $J$ or equivalently we view $M$ as an isotropic summand of $W$, $J/M$ as an isotropic summand of $M^\perp/M$, and $J^\perp/J$ as $(J/M)^\perp/(J/M)$.

We will check that they are equal by checking they are equal on the summand associated to $M \subset J$. 
If there is  no  $\sigma \in S_{2n}^\omega$ and $k \leq m \leq n$ with 
\begin{gather*} 
    J=\Span(\vec v_{\sigma(1)},\vec v_{\sigma(3)},\ldots,\vec v_{\sigma(2m-1)}) \cap W\text{ and }\\
    M=\Span(\vec v_{\sigma(1)},\vec  v_{\sigma(3)},\ldots,\vec v_{\sigma(2k-1)}) \cap W,
\end{gather*}
then both terms are zero. 
Otherwise, both terms are equal to

    $$ \chi(\sigma)[\vec v_{\sigma(1)},\vec v_{\sigma(3)}\ldots, \vec v_{\sigma(2k-1)} ] \otimes [\vec v_{\sigma(2k+1)},\vec v_{\sigma(2k+3)}\ldots, \vec v_{\sigma(2m-1)} ]  \otimes [\vec v_{\sigma(2m+1)},\vec v_{\sigma(2m+2)},\ldots, \vec v_{\sigma(2n)} ] . $$  
Thus $(\St^\omega,\Delta)$ is a coassociative comodule over $(\St,\Delta)$.
\end{proof}

\subsection{Duoidal comodule algebra structure} \label{compat}

Kupers--Rudenko--Sierra \cite[Section 2.3]{KRS1} checked that
 $(\St,\Delta,\mu)$ is a graded-commutative duoidal bialgebra in the duoidal category $(\Mod_{\VB(R)}(\bk),\otimes_\PP,\otimes_\LL)$ for $R$ a field. A similar argument works for $R$ a Dedekind domain.

\begin{theorem} \label{commutingSp}
 For $R$ a Dedekind domain, $(\St^\omega,\Delta,\mu)$ is a graded-commutative duoidal comodule algebra over $(\St,\Delta,\mu)$.
\end{theorem}

\begin{proof}
The only interesting condition of \autoref{duoidalcomodulealgebra} is \autoref{duoidalcomodulealgebra1}. To check that, we must show  that the following diagram commutes:  
\begin{center}   
        \begin{tikzcd}
        \St^\omega \otimes_\LL^\omega \St^\omega \arrow[d, "\mu"] \arrow[r, "\Delta \otimes_\LL^\omega \Delta"]&  (\St \otimes_\PP^\omega \St^\omega) \otimes_\LL^\omega (\St \otimes_\PP^\omega \St^\omega) \arrow[d, "(\mu \otimes_\PP^\omega \mu) \circ \xi"]  \\ 
            \St^\omega  \arrow[r, "\Delta"] & \St \otimes_\PP^\omega \St^\omega.
        \end{tikzcd}
    \end{center} 
Concretely, we must check that the following diagram commutes for every symplectic module $W$:

\begin{center}   
        \begin{tikzcd} 
         \bigoplus_{U \oplus V= W}   \St^\omega(U) \otimes \St^\omega(V) \arrow[d, "\mu"] \arrow[r, "\Delta \otimes_\LL^\omega \Delta"]& \bigoplus_{\substack{U \oplus V =W, \\ J \subset U, \, K \subset V}} \St(J) \otimes \St^\omega(J^\perp/J) \otimes \St(K) \otimes \St^\omega(K^\perp/K)   \arrow[d, "(\mu \otimes_\PP^\omega \mu) \circ \xi"]\\ 
            \St^\omega(W)  \arrow[r, "\Delta"] & \bigoplus_{Q \subset W} \St(Q) \otimes \St^\omega(Q^\perp/Q).
        \end{tikzcd}
    \end{center}
    Let $[\vec a_1,\vec b_1,\ldots,\vec a_k,\vec b_k] \otimes [\vec c_1,\vec d_1,\ldots,\vec c_j,\vec d_j]$ be an element of the top left. Traversing the diagram counter clockwise yields $$\Delta\Big(\mu\Big([\vec a_1,\vec b_1,\ldots,\vec a_k,\vec b_k] \otimes [\vec c_1,\vec d_1,\ldots,\vec c_j,\vec d_j]\Big)\Big)=\Delta\Big([\vec a_1,\vec b_1,\ldots,\vec a_k,\vec b_k,\vec c_1,\vec d_1,\ldots,\vec c_j,\vec d_j]\Big)$$
This is the signed sum of apartments associated to isotropic subsets of $$\{\vec a_1,\vec b_1,\ldots,\vec a_k,\vec b_k,\vec c_1,\vec d_1,\ldots,\vec c_j,\vec d_j\}$$ tensored with the symplectic apartment associated with the complementary symplectic basis. 

Traversing the diagram clockwise yields $$(\mu \otimes_\LL \mu) \big( \xi \big( \Delta([\vec a_1,\vec b_1,\ldots,\vec a_k,\vec b_k]) \otimes\Delta( [\vec c_1,\vec d_1,\ldots,\vec c_j,\vec d_j]) \big) \big). $$ Note that $\Delta([\vec a_1,\vec b_1,\ldots,\vec a_k,\vec b_k])$ is the signed sum of apartments associated to isotropic subsets of $$\{\vec a_1,\vec b_1,\ldots,\vec a_k,\vec b_k\}$$ tensored with the symplectic apartment associated with the complementary symplectic basis and similarly $\Delta( [\vec c_1,\vec d_1,\ldots,\vec c_j,\vec d_j])$ and has the analogous description. Applying the function $(\mu \diamond \mu) \circ \xi$ concatenates the symplectic apartments and the isotropic apartments. Since isotropic subsets of $$\{\vec a_1,\vec b_1,\ldots,\vec a_k,\vec b_k,\vec c_1,\vec d_1,\ldots,\vec c_j,\vec d_j\}$$ are the same as sets obtained by taking unions of isotropic subsets of $$\{\vec a_1,\vec b_1,\ldots,\vec a_k,\vec b_k \} \text{ and }\{\vec c_1,\vec d_1,\ldots,\vec c_j,\vec d_j\},$$ traversing the diagram clockwise and counterclockwise agree. 
\end{proof}

\subsection{Classical comodule algebra structure}

In this section, and whenever talking about $\HH$ and $\HA$, we will assume that $\bk$ is a field. This is done so that the K\"unneth map is an isomorphism.

\begin{definition}
    Let $R$ be a Dedekind domain. For each isomorphism class of finitely generated projective modules $[P]$, pick a representative. Let
    $$\HH_{[P],i}= H_i(\GL(P);\St(P) \otimes \bk),$$ 
    $$\HH_{n,i}= \bigoplus_{\rk([P])=n} \HH_{[P],i},$$ 
    $$\text{ and }\HH= \bigoplus_{n,i} \HH_{n,i},$$
   Let $$\HA_{n,i}=H_i(\Sp_{2n}(R);\St_{2n}^\omega(R) \otimes \bk) \text{ and }\HA=\bigoplus_{n,i} \HA_{n,i}.$$ 
\end{definition}

\subsubsection{Algebra structures on $\HH$ and $\HA$}

Let $P$ and $Q$ be finitely generated projective $R$-modules. The product $$\mu\colon  \St(P) \otimes \St(Q) \to \St(P \oplus Q)$$ induces a map $$\mu\colon  H_i(\GL(P);\St(P)) \otimes H_j(\GL(Q); \St(Q)) \to H_{i+j}(\GL(P \oplus Q); \St(P \oplus Q)).$$
By picking isomorphisms between $P$, $Q$, and $P \oplus Q$ with our chosen representatives of $[P]$, $[Q]$, and $[P \oplus Q]$, we obtain a product $$ \HH_{[P],i} \otimes \HH_{[Q],j} \to \HH_{[P \oplus Q],i+j}.$$ Note that this product does not depend on these choices of isomorphisms since conjugation by an element $g$ of a group $G$ induces the trivial automorphism on the group homology of $G$. In this way, $\mu$ induces a graded ring structure on $\HH$. Similarly, $\HA$ has a graded ring structure induced by $\mu\colon \St^\omega \otimes_\LL^\omega \St^\omega \to \St^\omega$. 

These products are graded-commutative in these sense that if $\alpha \in \HH_{n,i}$ and $\beta \in \HH_{m,j}$, then $$ \alpha \beta=(-1)^{nm+ij}\beta \alpha$$ an similarly for $\HA$.

\subsubsection{Coalgebra structure on $\HH$ and comodule structure on $\HA$}

The coproducts induce coalgebra and comodule structures on $\HH$ and $\HA$ respectively. Let $M$ be a finitely generated projective $R$-module. The coproduct $\Delta \colon \St \to \St \otimes_\PP \St $ induces a map $$\Delta\colon  H_i(\GL(M);\St(M))) \to H_i(\GL(M);\bigoplus_{M' \subseteq M } \St(M') \otimes \St(M/M') )).$$ We will simplify the codomain using Shapiro's lemma. Note that the orbit of $\GL(M)$ acting on the set of summands of $M$ is the set of isomorphism classes of $R$-modules of rank less than or equal to the rank of $M$ that are isomorphic to $M$ if they have the same rank as $M$. Thus, Shapiro's lemma provides a isomorphism \begin{multline*}H_i(\GL(M);\bigoplus_{M' \subseteq M } \St(M') \otimes \St(M/M') ))\\ \longrightarrow \bigoplus_{\substack{[M'],\, \rk(M')\leq \rk(M)\\ [M'] = [M] \text{ if } \rk(M') = \rk(M)}} H_i(\Stab_{M'}(\GL(M));\St(M') \otimes \St(M/M')).\end{multline*} 
The stabilizer of $M'$ in $\GL(M)$, is a parabolic subgroup with a natural split surjective map to the corresponding Levi subgroup:
$$\Stab_{\GL(M)}(M') \to \GL(M') \times \GL(M/M').$$ The action of the parabolic on $\St(M') \otimes \St(M/M')$ factors through the Levi so we obtain a map $$H_i(\GL(M);\St(M))) \to \bigoplus_{[M'],\, \rk(M')\leq \rk(M)} H_i( \GL(M') \times \GL(M/M') ;\St(M') \otimes \St(M/M')). $$
 After applying the inverse of the K\"unneth map, we obtain a coalgebra structure: $$\Delta \colon  \HH_{n,i} \to \bigoplus_{\overset{a+b=n}{k+j=i}} \HH_{a,j} \otimes \HH_{b,k}.$$ 

A similar discussion provides a comodule map $$\Delta \colon  \HA_{n,i} \to \bigoplus_{\overset{a+b=n}{k+j=i}} \HH_{a,j} \otimes \HA_{b,k}.$$ Here we use that if $J$ is an isotropic summand of a symplectic $R$-module $W$, then the Levi subgroup corresponding to parabolic subgroup $\Stab_{\Sp(W)}(J)$ is $$\GL(J) \times \Sp(J^\perp/J).$$ Again, the action of the parabolic on $ \St(J) \otimes \St^\omega(J^\perp/J)$ factors through the Levi.

\subsubsection{Hopf algebra structure on $\HH$ and comodule algebra structure on $\HA$}

Since $\St$ is a duoidal bialgebra, the product and coproduct on $\HH$ are compatible in the sense that they endow $\HH$ with the structure of a  bialgebra. Since the product on $\St$ is  graded-commutative, so is the product on $\HH$. Since $\HH_{0,0}$ is the $0$th homology of the trivial group with coefficients in $\bk=\St(0)$, $\HH_{0,0}=\bk$. Connected bialgebras are Hopf algebras (Milnor--Moore \cite[Proposition 8.2]{MM}). 

\begin{proposition}
    Let $\bk$ be a field and $R$ a Dedekind domain. Then $\HH$ is a connected  graded-commutative Hopf algebra.
\end{proposition}
Recall that a comodule algebra over a bialgebra $B$ is a monoid object in the category of $B$-comodules. The duoidal $\St$-comodule algebra structure on $\St^\omega$ gives the following which is a restatement of \autoref{main}.

\begin{theorem}
    Let $\bk$ be a field and $R$ a Dedekind domain. Then $\HA$ is a connected graded-commutative comodule algebra over the Hopf algebra $\HH$. 
\end{theorem}

\section{Freeness of comodule algebras}
\label{SecFree}
In this section, we prove a structure theorem for module coalgebras, which will give us a lower bound for the comodule algebra $\HA$ in the following section. 
Throughout this section, let $\bk$ be a field, $B$ be a graded connected bialgebra over $\bk$, and $C$ a graded connected $B$-module coalgebra. 

We denote by $W$ the image
\[W:= \im( B \cong B \otimes C_0 \hookrightarrow
B \otimes C \to C). \]
Note that $W$ is a graded $B$-submodule of $C$, because 
\[ b' \cdot ( b\cdot 1_C) = (b' \cdot b)  \cdot 1_C .\]
$W$ is also a graded subcoalgebra of $C$ because
\[ \Delta_C(b \cdot 1_C) = \Delta_B(b) \cdot \Delta_C(1_C) = \Delta_B(b) \cdot (1_C \otimes 1_C) \in W \otimes W.\]

Further, denote by $Q$ the graded quotient vector space
\[ Q:= C/W.\]
Let $q\colon C \to Q$ denote the canonical quotient map.
It follows that $Q$ is a graded $B$-module. Furthermore, $Q$ inherits a $C$-comodule structure, as
\[ \Delta_C( W) \subseteq W \otimes W \subseteq C \otimes W\]
and thus
\begin{align*} \Delta_Q\colon Q &\longrightarrow C \otimes Q  \\
c + W &\longmapsto (\id \otimes q) \Delta_C(c)
\end{align*}
is well defined. Thus $Q$ is a $B$-module $C$-comodule.

The $C$-comodule primitives of $Q$ are defined as
\[P(Q) := \{ c+W\in Q \mid \Delta_Q(c+W) = 1_C\otimes (c+W)\}.\]

\begin{lemma}
    $W$ is a subspace of
    \[ \{ c\in C \mid \Delta_C(c) - 1_C \otimes c \in C \otimes W\}\]
    and 
    \[ P(Q) = \{ c\in C \mid \Delta_C(c) - 1_C \otimes c \in C \otimes W\}/W.\]
\end{lemma}

\begin{proof}
    Since $\Delta_C(w) \in W \otimes W$ for all $w$, so is $\Delta_C(w) - 1_C \otimes w$.

    It remains to prove that $\Delta_C(c) -1_C \otimes c \in C \otimes W$ if and only if $\Delta_Q(c+W) = 1_C \otimes (c+W)$.

    This can be seen because 
    \[ \Delta_Q(c+W) - 1_C \otimes (c + W) = (\id \otimes q)(\Delta_C(c) - 1_C \otimes c) \]
    and $(\id \otimes q)(x) =0 $ if and only if $x\in C \otimes W$.
\end{proof}

We will now prove our main structure theorem.

\begin{theorem} \label{FreenessModulecoAlgebra}
    There is a graded injection from $W \otimes (P(Q) \oplus C_0)$ into $C$.
\end{theorem}

The main part of the proof is in the following lemma.

\begin{lemma}
    Let $s\colon W \to B$ be a graded section of the surjection $B \to W$. Then
    \[ \Phi\colon W \otimes P(Q) \xrightarrow{s \otimes \operatorname{incl}} B \otimes Q\stackrel{\mu_Q}{\longrightarrow} Q \]
    is a graded injection.
\end{lemma}

\begin{proof}
    It is easy to observe that $\Phi$ is graded as all maps involved are graded.

    Assume that $x\in W \otimes P(Q)$ is nonzero. We need to show that $\Phi(x) \neq 0$.
    Let us assume that 
    \[ x = x_0 + \dots +x_n \]
    with $x_i \in W_i \otimes P(Q)$ and $x_n \neq 0$.

    We calculate that for any $w_i \otimes (c+W) \in W_i \otimes P(Q)$ we get
    \begin{align*}
        &\Delta_Q(\Phi(w_i \otimes (c+W)))\\
        =&\Delta_Q((s(w_i) \cdot (c+W))) \\
        =& \Delta_B(s(w_i))\cdot \Delta_Q(c+W) \\
        =& \Delta_B(s(w_i))\cdot (1_C \otimes (c+W)) \\
        =& ( s(w_i) \otimes 1_B + v) \cdot (1_C \otimes (c+W)) \\
        =& w_i \otimes (c+W) + u, 
    \end{align*}
    where $v \in B_{<i} \otimes B$ with $B_{<i} = B_0 \oplus \dots \oplus B_{i-1}$ and $u \in C_{<i}\otimes Q$ with $C_{<i} = C_0 \oplus \dots \oplus C_{i-1}$. By linearity, we conclude that
    \[\Delta_Q(\Phi(x)) = x_n + u',\]
    where $u'\in C_{<n}\otimes Q$.   Since $x_n+ C_{<n} \otimes Q$ is nonzero, $\Phi(x) \neq 0$.
\end{proof}

\begin{proof}[Proof of \autoref{FreenessModulecoAlgebra}]
    Let $r\colon Q \to C$ be a graded section of $q\colon C \to Q$. Then $W \oplus r(Q)=C$. Thus
    \begin{align*}
        \Psi\colon W \otimes (P(Q) \oplus C_0) & \longrightarrow C\\
        w \otimes (c+W) &\longmapsto r(s(w)\cdot c + W)\\
        w \otimes 1_C &\longmapsto w
    \end{align*}
    is a monomorphism.
\end{proof}

For more applicability, we finish this section by identifying a particular subspace of $P(Q)$.

\begin{lemma}
    $(P(C)+W)/ W$ is a subspace of $P(Q)$, where $P(C)$ are the primitives of $C$.
\end{lemma}

\begin{proof}
    This follows immediately because $\Delta_C(c) - 1_C \otimes c = c \otimes 1_C \in C \otimes W$ for every $c\in P(C)$.
\end{proof}

The last lemma has the following direct consequence.

\begin{corollary}
    $W \otimes [C_0 \oplus P(C)/(P(C) \cap W)]$ injects into $C$.
\end{corollary}

In our applications, we actually want to use the dual statement for $B$-comodule algebras $A$. To this end, the dual of $W$ is
\[ V := \im( A \stackrel{\Delta_A}{\longrightarrow} B \otimes A \stackrel{\id \otimes \epsilon_A}{\longrightarrow} B \otimes A_0 \cong B)\]
and the dual of $P(C) = Q(A) = A_+/A^2_+$. We denote composition $$ \id \otimes \epsilon_A \circ  \Delta_A$$ by $\delta:A\to B$ and call this the structure homomorphism of the $B$-comodule algebra structure of $A$. This leads to the following statement for a graded connected bialgebra $B$ and a graded connected $B$-comodule algebra $A$ both of finite type. 

\begin{corollary} \label{ComoduleAlgFree}
    $A$ surjects onto $V \otimes (A_0 \oplus \ker(\delta \colon Q(A) \to Q(V))$.
\end{corollary}

\section{Consequences for cohomology of symplectic groups} \label{ApSm}

The goal of this section is to establish lower bounds for the high-dimensional cohomology of symplectic groups of number rings.

\subsection{A source of  indecomposables}

Borel--Serre proved symplectic groups and general linear groups are virtual duality groups.

\begin{theorem}[Borel--Serre] \label{BSduality} Let $R$ be a number ring with field of fractions $F$. Let $r_1$ be the number real embeddings of $F$ and $r_2$ pairs of complex embeddings of $F$. Let $M$ be a projective $R$-module of rank $n$ and let $W$ be a symplectic module of genus $n$. The virtual cohomological dimension of $\GL(M)$ is $$\nu^{\GL}_n :=r_1 {{n+1} \choose 2} +r_2 n^2 -n$$ and the virtual cohomological dimension of $\Sp(W)$ is $$\nu_{2n}^\omega:= r_1 (n(n+1))+r_2(n(2n+1))-n.$$ The virtual dualizing module of $\Sp(W)$ is $\St^\omega(W)$ and $\St(M)$ is the dualizing module of $\SL(M)$. A tensor product of $\St(M)$ with a one-dimensional representation is the dualizing module of $\GL(M)$; see 
  Putman--Studenmund \cite[Theorem C]{PutStu}.
\end{theorem}
Recall 
$$\HH :=\bigoplus_{[P],i} H_i(\GL(P);\St(P \otimes_R F))$$ and 
$$\HA :=\bigoplus_{n,i } H_i(\Sp_{2n}(R);\St^\omega(F^{2n})).$$ 
 Let $\delta\colon\HA \to \HH$ denote the structure homomorphism 
 and let  $\VV=\im(\delta) \subseteq \HH$. In particular, $\delta$ is given by the composition 
 \begin{multline*} H_i(\Sp_{2n}(R); \St^\omega(F^{2n})) \\\longrightarrow H_i(\Sp_{2n}(R); \bigoplus_{L\in \fL(R^{2n})} \St(L)) \cong \bigoplus_{[L], L\in \fL(R^{2n})} H_i(\Stab_{\Sp_{2n}(R)}(L);\St(L \otimes_R F)) \\\longrightarrow \bigoplus_{[L], \rk(L) = n}H_i(\GL(L);\St(L \otimes_R F)) , \end{multline*}
where $\fL(W)$ denotes the set of Lagrangians of a symplectic $R$-module $W$, the first map is induced by the coproduct $\bigoplus_{L \in \fL(R^{2n})}\Delta_{L\subset R^{2n}}$, and the second map by the projection $\Stab_{\Sp_{2n}(R)}(L) \to \GL(L)$.

We will prove that the kernel $$\ker(\delta\colon Q(\HA) \to Q(\VV))$$ is large for degree reasons. 

\begin{lemma} \label{sourceofindecomposables} Let $\bk$ be a field of characteristic zero. Let $R$ be a number ring with field of fractions $F$. Let $r_1$ be the number of real embeddings of $F$ and $r_2$ be the number of pairs of complex embeddings of $F$. Let 
$$S=\{(n,i) \,| \, i= \nu_{2n}^{\omega} -j \text{ with } j \leq 2(r_1+2r_2)(n-1)-1 \} \cup \{(1,i) \, | \, i > r_1+r_2-1  \} \subset \N^2.$$ 
    There is a natural injection $$\bigoplus_{(n,i) \in S} \HA_{n,i} \hookrightarrow  \ker(\delta\colon Q(\HA) \to Q(\VV)).$$ 
\end{lemma}

\begin{proof}
    Suppose $\alpha \neq 0 \in \HA_{n,i}$ is decomposable. Then $n>1$ and $\alpha$ is a sum of products of elements $\beta$ and $\gamma$ of the form $\beta \neq 0 \in \HA_{a,b}$ and a $\gamma \neq 0 \in \HA_{c,d}$ with $n=a+c$, $i=b+d$, and $a,c>0$. By \autoref{BSduality}, $b \leq \nu_{2a}^\omega$ and $d \leq \nu_{2c}^\omega$. Thus, $$i \leq \nu_{2a}^\omega+\nu_{2c}^\omega \leq \nu_{2}^\omega+\nu_{2n-2}^\omega = \nu_{2n}^\omega - 2(r_1+2r_2)(n-1).$$ The second inequality follows from convexity of the function $n \mapsto \nu^\omega_{2n}.$    
    Hence
    $$\bigoplus_{(n,i) \in S} \HA_{n,i} \hookrightarrow  Q(\HA).$$ \autoref{BSduality} implies $\HH_{n,i}=0$ for $i>\nu_n$. In particular, $\HH_{n,i} = 0$ if $i > \nu_{2n}^\omega - 2(r_1+2r_2)(n-1) \ge \nu_n$ and if $n=1$ and $i> \nu_1 =  r_1+r_2-1$. Thus $\HH_{n,i} = 0$ if $(n,i)\in S$. As $Q(\VV)$ is a quotient of $\HH$, it also vanishes in this range. The claim follows. 
\end{proof}

The best known homological stability range for symplectic groups is due to Sierra–Wahl \cite{SierraWahl}. 

\begin{theorem}[Sierra–Wahl] \label{SW}
    Let $R$ be a Dedekind domain and $i \leq \frac{2n-8}{3}$. Then $H^i(\Sp_{2n}(R)) \cong H^i(\Sp_{2n+2}(R))$. 
\end{theorem}

Borel \cite{Bor} computed the stable cohomology. 

\begin{theorem}[Borel] \label{Borel}  Let $R$ be a number ring with field of fractions $F$. Let $\bk$ be a field of characteristic zero. Let $r_1$ be the number of real embeddings of $F$ and $r_2$ the number of pairs of complex embeddings of $F$.
Then the stable cohomology $H^*(\Sp_{\infty}(R))$ is the graded algebra $$\bk[\{x_{i,j} \mid i\ge 1, 1\le j \le r_1\}] \otimes {\bigwedge}_{\bk}[\{y_{a,b}\mid a\ge 1, 1\le b \le r_2 \}]$$ with $|x_{i,j}|=4i-2$ and $|y_{a,b}|=4a-1$.
\end{theorem}

\subsection{Lower bounds for the image of the structure homomorphism in the vcd}

Given a symplectic module $W$, recall that $\fL(W)$ denotes the set of Lagrangians of $W$. The structure homomorphism $\delta:\HA \to \HH$ is induced by a map
$$\delta^\omega\colon \St^\omega(W) \longrightarrow \bigoplus_{L \in \fL(W)} \St(L)$$ defined as the sum of the maps $\Delta_{L \subset W}$.

For $M$ a rank $d$-module, let $$\delta^{d-1} \colon \St(M) \longrightarrow \bigoplus_{\substack{N \subseteq M\\\rk(N)=d-1}} \St(N) \otimes \St(M/N) $$ be the sum of the maps $\Delta_{N \subseteq M}$ with $\rk(N)=d-1$. Let
\begin{multline*}D^{d-1} \colon \bigoplus_{\substack{J_d \subset J_{d+1} \subset \dots \subset J_n\\\rk(J_i) =i\\J_n \in \mathfrak L(W)}} \St(J_d) \otimes \St(J_{d+1}/J_d) \otimes \dots \otimes \St(J_n/J_{n-1}) \\
\longrightarrow\bigoplus_{\substack{J_{d-1} \subset J_d \subset \dots \subset J_n\\\rk(J_i) =i\\J_n \in \mathfrak L(W)}} \St(J_{d-1}) \otimes \St(J_d/J_{d-1}) \otimes \dots \otimes \St(J_n/J_{n-1}) 
\end{multline*}
be the map induced by $\delta^{d-1}$. Let 
$$\delta^{\mathrm{total}}\colon \St^\omega(W) \longrightarrow \bigoplus_{\substack{J_1 \subset J_2 \subset \dots \subset J_n\\\rk(J_i) =i\\J_n \in \mathfrak L(W)}} \St(J_1) \otimes \St(J_2/J_1) \otimes \dots \otimes \St(J_n/J_{n-1})$$ be the composition $$\delta^{\mathrm{total}}:=D^1 \circ D^2 \circ \dots\circ D^{n-1} \circ   \delta^\omega.$$ Since $J_1, J_2/J_1,\ldots,J_n/J_{n-1}$ are rank one, their Steinberg modules are naturally isomorphic to $\bk$. The indexing set of the direct sum is the set of full flags of isotropic summands of $W$. Thus,  we may view $\delta^{\mathrm{total}}$ as a map to the top cellular chains on the symplectic Tits building: $$\delta^{\mathrm{total}}\colon \St^\omega(W) \to \widetilde C_{n-1}(\cT^\omega(W)).$$ Here $n$ is the genus of $W$.

\begin{lemma}\label{lemma:deltatotal}
    The map $\delta^{\mathrm{total}}\colon \St^\omega(W) \to \widetilde C_{n-1}(\cT^\omega(W))$ agrees with the inclusion $$\widetilde H_{n-1}(\cT^\omega(W)) \to \widetilde C_{n-1}(\cT^\omega(W)).$$
\end{lemma}
\begin{proof}
Let $(\vec v_1,\ldots,\vec v_{2n})$ be an ordered  symplectic basis of $W \otimes_R F$. After unpacking the definition, we see that $\delta^{\mathrm{total}}([\vec v_1,\ldots,\vec v_{2n}])=$ $$\bigoplus_{\sigma \in S^\omega_{2n}} \chi(\sigma) [\vec v_{\sigma(1)}] \otimes [\vec v_{\sigma(3)}] \otimes \ldots \otimes [\vec v_{\sigma(2n-1)}]. $$  Here $[\vec v_{\sigma(1)}] \otimes [\vec v_{\sigma(3)}] \otimes \ldots \otimes [\vec v_{\sigma(2n-1)}]$ is in the summand associated to the flag $$\Span(\vec v_{\sigma(1)}) \subset \Span(\vec v_{\sigma(1)},\vec v_{\sigma(3)} ) \subset \ldots \subset \Span(\vec v_{\sigma(1)},\ldots,\vec v_{\sigma(2n-1)}). $$ The signed sum of these flags is exactly a symplectic apartment.
\end{proof}

Let $W$ have genus $n$. The quotient map $\cT^\omega(W) \to \cT^\omega(W)/\Sp(W)$ is invariant and hence induces a map $$H_0(\Sp(W);\St^\omega(W)  ) \to \widetilde H_{n-1}( \cT^\omega(W)/\Sp(W)). $$ We now recall the following result of Br\"uck--Himes \cite[Theorem 1.2 and Section 3]{BH}.

\begin{theorem}[Br\"uck--Himes]\label{bruckhimes} Let $R$ be a Dedekind domain with class number $c$ and let $W$ be a symplectic module with genus $n$. There is a map $$\pi \colon \widetilde H_{n-1}( \cT^\omega(W)/\Sp(W)) \to \bk^{(c-1)^n}$$ such that the composition
$$  H_0(\Sp(W);\St^\omega(W)  )  \to \widetilde H_{n-1}( \cT^\omega(W)/\Sp(W)) \to \bk^{(c-1)^n}$$ is surjective.
\end{theorem}

We obtain the following corollary on the image of the structure homomorphism. 

\begin{corollary} \label{Cor2} Let $R$ be a Dedekind domain with class number $c$. Let $\bk$ be a field. There is a surjection $$\VV_{n,0} \twoheadrightarrow \bk^{(c-1)^n}.$$

\end{corollary}

\begin{proof} 
   \autoref{lemma:deltatotal} implies that  
    \[\St^\omega(R^{2n})  \stackrel{\delta^{\mathrm{total}}}{\longrightarrow} C_{n-1}(\cT^\omega(R^{2n})) \hooklongrightarrow \widetilde C_*(\cT^\omega(R^{2n}))\] 
    is a (degree $n-1$) quasiisomorphism and thus induces an isomorphism
    \[ H_0(\Sp_{2n}(R);\St^\omega(F^{2n})) \longrightarrow H_0(\Sp_{2n}(R);\widetilde C_*(\cT^\omega(F^{2n})).\]
    There is an equivariant map 
    \[\widetilde C_*(\cT^\omega(R^{2n})) \longrightarrow \widetilde C_*(\cT^\omega(R^{2n})/\Sp_{2n}(R))\]
    and $\widetilde C_*(\cT^\omega(R^{2n})/\Sp_{2n}(R))$ is quasiisomorphic to $\widetilde H_{n-1}(\cT^\omega(R^{2n})/\Sp_{2n}(R))$ on which $\Sp_{2n}(R)$ acts trivially. \autoref{bruckhimes} implies that the composition
    \[ H_0(\Sp_{2n}(R);\St^\omega(F^{2n})) \longrightarrow H_0(\Sp_{2n}(R);\widetilde H_{n-1}(\cT^\omega(R^{2n})/\Sp_{2n}(R))) = \widetilde H_{n-1}(\cT^\omega(R^{2n})/\Sp_{2n}(R)) \stackrel{\pi}{\longrightarrow}\bk^{(c-1)^n} \]
    is a surjection.

    Observe that the map
    \[ H_0(\Sp_{2n}(R); \bigoplus_{L\in \fL(R^{2n})} \St(L)) \cong \bigoplus_{[L], L\in \fL(R^{2n})} H_0(\Stab_{\Sp_{2n}(R)}(L);\St(L \otimes F)) \longrightarrow \bigoplus_{[L], \rk(L) = n}H_0(\GL(L));\St(L \otimes F)) \]
    is an isomorphism because $\Stab_{\Sp_{2n}(R)}(L)$ acts on $\St(L)$ through the projection $\Stab_{\Sp_{2n}(R)}(L) \to \GL(L)$.

    Thus because $\delta^{\mathrm{total}}$ factors through
    \[ \delta^\omega\colon \St^\omega(F^{2n}) \longrightarrow \bigoplus_{L\in \fL(R^{2n})} \St(L),\]
    the surjective composition
    \[ \HA_{n,0} \cong H_0(\Sp_{2n}(R);\St^\omega(F^{2n})) \longrightarrow \bk^{(c-1)^n} \]
    factors through $\delta\colon \HA_{n,0} \to\HH_{n,0}$ and thus through $\VV_{n,0}$.
\end{proof}
 
\subsection{Lower bounds parallel to the vcd}

Using items (1), (2), (3) of the following theorem, we obtain \autoref{corClassNumber} except in the case $m=0$ which is Br\"uck--Himes \cite[Theorem 1.1]{BH}. 

\begin{theorem} \label{corClassNumberBetter}
Let $\bk$ be a field of characteristic zero. Let $R$ be a number ring with field of fractions $F$ and class number $c$. Let $r_1$ be the number real embeddings of $F$ and $r_2$ pairs of complex embeddings of $F$.
Let $$S=\{(n,i) \,| \, i= \nu_{2n}^{\omega} -j \text{ with } j \leq 2(r_1+2r_2)(n-1)-1 \} \cup \{(1,i) \, | \, i > r_1+r_2-1  \} \subset \N^2.$$ 
\begin{enumerate}
    \item 
  There is a surjection from $\HA$ onto  $$ \Big( \bigoplus_n \VV_{n,0} \Big) \otimes \Big( \bigoplus_{(n,i) \in S} \HA_{n,i} \Big) $$ 

  \item  $\dim_\bk (\VV_{n,0}) \geq (c-1)^n$,
    \item (Borel--Serre) \, $\HA_{n,\nu^\omega_{2n}} \cong \bk$,
  \item (Borel--Serre, Sierra--Wahl) \, For $i = \nu_{2n}^\omega - j $ with $j \leq \frac{2n-8}{3}$, $\HA_{n,i} \cong H^j(\Sp_\infty(R ))$,
  \item (Borel) \, $H^*(\Sp_\infty(R)) \cong  \bk[\{x_{i,j} \mid i\ge 1, 1\le j \le r_1\}] \otimes {\bigwedge}_{\bk}[\{y_{a,b}\mid a\ge 1, 1\le b \le r_2 \}]$ with $|x_{i,j}|=4i-2$ and $|y_{a,b}|=4a-1$.
\end{enumerate}

\end{theorem}

\begin{proof}
\begin{enumerate}

    \item This combines \autoref{ComoduleAlgFree} and \autoref{sourceofindecomposables}.
    \item This is \autoref{Cor2}.
    \item Borel--Serre duality \autoref{BSduality} implies  that $$H_{\nu^\omega_{2n}}(\Sp_{2n}(R);\St^\omega(F^{2n})) \cong H^0(\Sp_{2n}(R)) \cong \bk.$$
    \item This combines Borel--Serre duality \autoref{BSduality} and Sierra–Wahl's homological stability result \autoref{SW}. 
    \item This is Borel's calculation of the stable rational cohomology of symplectic groups \autoref{Borel}.
\end{enumerate}
\end{proof}

Items (4) and (5) can be used to obtain more information.  

\begin{example} \label{codim96}
     If $R$ is the ring of integers in $\Q[\sqrt{-23}]$, then $r_1=0$, $r_2=1$, and $c=3$, $$H^3(\Sp_{18}(R);\bQ) \cong H^3(\Sp_{\infty}(R);\bQ) \cong \Q,$$ and $\nu^\omega_{18}=162$. We conclude that  
$$\dim_\Q H^{2n^2-159}(\Sp_{2n}(R);\Q) \geq 2^{n-9}$$ for $n \geq 9$.
\end{example}

\begin{remark}
 The proof of  \autoref{Cor2} implies that there is a surjection $\VV_{n,0} \to \HA_{n,0}$. It would be interesting to investigate if one could use this observation and ideas from \cite{MPWY} to construct lower bounds for the rational cohomology of symplectic groups of number rings with class number equal to $1$ that are not Euclidean.
\end{remark}

\begin{remark}\label{FrankeRem}
    Using Franke \cite[Section 7]{Franke} as in Speh--Venkataramana \cite[Section 4]{SpehVenk} and  Giansiracusa--Kupers--Tshishiku \cite[Theorem 13]{BenaSanderJeff}, one can describe the image of $H^*(\Sp_\infty (R))$ in $H^*(\Sp_{2n}(R))$. Let us define 
    \[ A_n = \bk[c_1,\dots,c_n]/I\]
    where $I$ is the ideal generated by the polynomials
    \[ \sum_{a+b = k} (-1)^b c_ac_b -\delta_{k,0}\]
    for all $k$ with the convention that $c_0=1$ and $c_i=0$ if $i \not\in[0,n].$
    The ring map $\bk[x_1,x_2,\dots]\to A_n$ sending $x_i$ to $c_{2i-1}$ is surjective. Identify $(A_n)^{\otimes r_1}$ with the corresponding quotient of 
    \[\bk[\{c_{i,j} \mid 1\le i \le n, 1 \le j \le r_1\}],\]
    then 
    \[\bk[\{x_{i,j} \mid 1\le i, 1\le j \le r_1\}]  \longrightarrow (A_n)^{\otimes r_1}\]
    given by $x_{i,j}$ sent to $c_{2i-1,j}$ is surjective and
    the image of $H^*(\Sp_\infty (R))$ in $H^*(\Sp_{2n}(R))$ is isomorphic to
    \[ (A_n)^{\otimes r_1} \otimes {\bigwedge}_{\bk}[\{y_{a,b}\mid 1\le a\le n, 1\le b \le r_2 \}]\bigg)/ (c_{n,1} \cdots c_{n,r_1} \otimes y_{n,1} \wedge \dots \wedge y_{n,r_2}).\]
    And one can find a basis given by 
    \[ \prod c_{i,j}^{\alpha_{i,j}} \otimes \prod y_{a,b}^{\beta_{a,b}}\]
    with $\alpha_{i,j},\beta_{a,b} \in \{0,1\}$ such that the product $\alpha_{n,1}\dots \alpha_{n,r_1}\beta_{n,1} \dots\beta_{n,r_2} = 0$.

    This gives a lower bound for $H^*(\Sp_{2n}(R))$, many of which are in the range 
    $$*\ge \nu^\omega_{2n}-  r_1(2n-3) - r_2(4n-5)$$
    but outside the stable range and thus improves the lower bounds established in \autoref{corClassNumberBetter} (4) and (5).

    Furthermore, one can adapt \cite[Proof of Theorem 1.2 in Section 4.2]{BCGP} to show that all these images of stable classes are indecomposable, which would improve \autoref{sourceofindecomposables}.
\end{remark}

\section{Consequences for cohomology of congruence subgroups}
\label{ApCong}
In this section, we establish lower bounds for the cohomology of congruence subgroups of $\Sp_{2n}(\Z)$.

\begin{definition}
    Let $W$ be a symplectic $R$-module and $I$ an ideal of $R$. Let $\Sp(W,I)$ denote the kernel of $$\Sp(W) \to \Sp(W/I).$$ For $W=R^{2n}$, we denote $\Sp(R^{2n},I)$ by $\Sp_{2n}(I)$. Similarly defined congruence subgroups of general linear groups. 
\end{definition}

\begin{definition} Let $U$ be an $m$-dimensional vector space over a field $F$. A $\pm$-orientation on $U$ is a choice of basis of the one-dimensional $\F$-vector space $\wedge^m_F U$ up to multiplication by $\pm 1$. 
\end{definition}

If $(\vec b_1,\ldots,\vec b_n)$ is an ordered basis of a free $\Z$-module $N$, then $\vec b_1 \wedge \ldots \wedge \vec b_n$ gives a $\pm$-orientation on the $\F_p$-vector space $N/p$. If $T\colon N\to N$ is a $\Z$-linear isomorphism, then $$\vec b_1 \wedge \ldots \wedge \vec b_n\text{ and }T(\vec b_1) \wedge \ldots \wedge T(\vec b_n)$$ differ by the ratio $\det(T)= \pm 1$. Thus, $\beta$ and $\beta'$ induce the same $\pm$-orientation on $N/p$. In this way, $N/p$ has a natural $\pm$-orientation induced by the $\Z$-module $N$.

The following is a corollary of 
Capovilla-Searle \cite[Proposition 4.11]{Fabio}.

\begin{theorem}[Capovilla-Searle] \label{FabioQuotient}
    Let $W$ be a symplectic $\Z$-module. Let $J_1$ and $J_2$ be isotropic summands of $W$. Then $J_1$ and $J_2$ are in the same orbit of $\Sp(pW)$ if and only if $J_1/p=J_2/p$ and the induced $\pm$-orientations agree. Additionally, every $\pm$-oriented isotropic subspace of $W/p$ lifts to an isotropic subspace of $W$.
\end{theorem}

The following is Borel--Serre dual to \autoref{CorCong}. Note that since $\Sp_{2n}(p\Z)$ is torsion free for $p \geq 3$,  Borel--Serre duality applies without restrictions on $\bk$.

\begin{theorem}
    Let $\bk$ be a field and $p$ be an odd prime. Then 
$$\dim_\bk H_i(\Sp_{2n+2}(p\Z);\St_{2n+2}^\omega(\Q)) \geq  p^{2n+1} \left(\frac{p-1}{2}\right) \dim_\bk H_i(\Sp_{2n}(p\Z);\St_{2n}^\omega(\Q)).$$ 
\end{theorem}

\begin{proof}

Let $\cS$ denote the set of genus-$1$ symplectic summands of $\Z^{2n+2}$. Let $\cJ$ denote the set of rank-$1$ isotropic summands of $\Z^{2n+2}$. Consider the composition of the following restrictions of the product and the coproduct 
$$\bigoplus_{U \in \cS} \St^\omega(U) \otimes \St^\omega(U^\perp) \to \St^\omega(\Z^{2n+2}) \to \bigoplus_{J \in \cJ} \St(J) \otimes \St^\omega(J^\perp/J). $$ 
Applying the functor $H_i(\Sp_{2n+2}(p\Z) ; - )$ and Shapiro's lemma gives maps 
\begin{align*}
&\bigoplus_{[U] \in \cS/\Sp_{2n+2}(p\Z)} H_i(\Stab_{\Sp_{2n+2}(p\Z)}(U) ; \St^\omega(U) \otimes  \St^\omega(U^\perp)) \\
\longrightarrow &\phantom{\bigoplus_{[U] \in \cS/\Sp_{2n+2}(p\Z)}}H_i(\Sp_{2n+2}(p\Z); \St^\omega(\Z^{2n+2})) \\
\longrightarrow&\bigoplus_{[J] \in \cJ/\Sp_{2n+2}(p\Z)} H_i(\Stab_{\Sp_{2n+2}(p\Z)}(J) ; \St(J) \otimes \St^\omega(J^\perp/J)).
\end{align*}
There is a natural isomorphism:
$$\Stab_{\Sp_{2n+2}(p\Z)}(U) \cong \Sp(U,p\Z) \times \Sp(U^\perp,p\Z).$$ The projection from the parabolic to the Levi gives a natural surjective map: $$ \Stab_{\Sp_{2n+2}(p\Z)}(J) \to \GL(J,p\Z) \times \Sp(J^\perp/J,p\Z). $$

After applying these maps and the Künneth isomorphism, we obtain maps 
\begin{align*}
&\bigoplus_{[U] \in \cS/\Sp_{2n+2}(p\Z)} H_0(\Sp( U,p\Z) ; \St^\omega(U)) \otimes H_i(\Sp( U^\perp,p\Z ); \St^\omega(U^\perp))\\
\longrightarrow&\phantom{\bigoplus_{[U] \in \cS/\Sp_{2n+2}(p\Z)} }H_i(\Sp_{2n+2}(p\Z);\St^\omega(\Z^{2n+2}))  \\
\longrightarrow &\bigoplus_{[J] \in \cJ/\Sp_{2n+2}(p\Z)} H_0(\GL(J,p\Z) ; \St(J)) \otimes H_i(\Sp(J^\perp/J,p\Z); \St^\omega(J^\perp/J)). 
\end{align*}
Call the composition of the above two maps $\Phi$ and let $\Phi_{[J]}$ be the projection onto the summand associated to $[J] \in \cJ/\Sp_{2n+2}(p\Z)$. 

Fix $U \in \cS$ and consider $J_1, J_2 \subset U$ rank-$1$ isotropic summands of $U$ that are not in the same orbit of $\Sp_{2n+2}(p\Z)$ and let $J$ be another rank-$1$ isotropic summand of $U$. For $i=1,2$, let $\vec j_i \in J_i \otimes \Q$ be a basis of $J_i$. Let $\alpha \in H_i(\Sp(U^\perp,p\Z) ; \St^\omega(U^\perp)) $.
\autoref{commutingSp} implies 
\begin{equation}
\begin{split}
\Phi_{J}([\vec j_1,\vec j_2] \otimes \alpha) =[\vec j_1]\otimes \alpha \text{ if } [J]=[J_1] \\
\Phi_{J}([\vec j_1,\vec j_2] \otimes \alpha) = -[\vec j_2] \otimes \alpha  \text{ if } [J]=[J_2] \\
\Phi_{J}([\vec j_1,\vec j_2] \otimes \alpha) = 0  \text{ otherwise.}
\end{split}
\label{PhiJ}
\end{equation}
Here we use the natural isomorphisms $$J^\perp/J \cong U^\perp $$ to identify the homology groups $$H_i(\Sp( J^\perp/J,p\Z);\St^\omega(J^\perp/J))  \cong H_i(\Sp( U^\perp,p\Z);\St^\omega(U^\perp))$$ and will view $\alpha$ as an element of either group.

Let $\overline v_1,\ldots, \overline v_{2n+2}$ be the standard symplectic basis of $\F_{p}^{2n+2}$. 
Let $\cJ'$ be the set of isotropic summands $J$ of $\Z^{2n+2}$ with $$J/p =\Span(\overline v_2 + \sum_{k \neq 2} x_k \overline v_k).$$ 
By \autoref{FabioQuotient}, there are $p^{2n+1}(p-1)/2$ many elements in the image of $\cJ'$ in $\cJ/\Sp_{2n+2}(p\Z)$. The  factor of $p^{2n+1}$ comes from the choices of $x_i \in \F_p$ and the factor of $(p-1)/2$ comes from the choice of orientation. 

Let $\Phi'$ be postcomposition of $\Phi$ with the projection onto the summands coming from elements of $\cJ/\Sp_{2n+2}(p\Z)$ in the image of $\cJ'$. The codomain of $\Phi'$ is
$$\bigoplus_{[J] \in \cJ/\Sp_{2n+2}(p\Z)\text{ with } J \in \cJ'} H_0(\GL(J,p\Z) ; \St(J)) \otimes H_i(\Sp(J^\perp/J,p\Z); \St^\omega(J^\perp/J))$$ and
has dimension 
$$p^{2n+1}\left(\frac{p-1}{2}\right) \dim_\bk H^{\nu^\omega_{2n}-i}(\Sp_{2n}(p\Z)).$$ Since $\Phi'$ factors through $H^{\nu^\omega_{2n+2}-i}(\Sp_{2n+2}(p\Z))$, it suffices to show $\Phi'$ is surjective to establish our lower bounds for the dimension of $H^{\nu^\omega_{2n+2}-i}(\Sp_{2n+2}(p\Z))$. 

We now prove that $\Phi'$ is surjective. Fix $J \in \cJ'$ and $\alpha \in H_i(\Sp(J^\perp/J,p\Z); \St^\omega(J^\perp/J))$. Pick $\vec j \in \Q^{2n+2}$ a vector with $\Span_\Q(\vec j)=J \otimes \Q$. View $[\vec j]$ as the generator of $\St(J)$. Elements of the form $$[\vec j] \otimes \alpha \in H_0(GL(J,p\Z);\St(J)) \otimes H_i(\Sp(J^\perp/J,p\Z);\St^\omega(J^\perp/J))$$ generate the codomain of $\Phi'$.

Using Capovilla-Searle \cite[Lemma 4.5]{Fabio}, we can find a symplectic summand $U \subseteq \Z^{2n+2}$  containing $J$ with $$U/p=\Span_{\F_p}(\overline v_1) \oplus J/p.$$ 

Let $\vec k \in U \otimes \Q$ be a vector with $$\Big(\Span_\Q(\vec k) \cap \Z^{2n+2}\Big)/p=\Span(\overline v_1).$$ Since $U^\perp \cong J^\perp/J$, we may view $\alpha$ as an element of $H_i(\Sp(U^\perp,p\Z); \St^\omega(U^\perp))$. Since $U \in \cS $ and $\Span_\Q(\vec k) \cap \Z^{2n+2} \notin \cJ'$, \autoref{PhiJ} implies 
$$\Phi'([\vec j,\vec k] \otimes \alpha)=[\vec j]\otimes \alpha.$$ Thus $\Phi'$ is surjective. 
\end{proof}




\bibliographystyle{amsalpha}
\bibliography{ref}

\end{document}